\documentclass[letter,twocolumn]{autart}
\usepackage{natbib}

\usepackage{amsmath,mathrsfs}
\usepackage{amssymb}
\usepackage{color}
\usepackage[all]{xy}
\usepackage{graphicx,xspace,bm}
\usepackage[caption = false]{subfig}
\usepackage{url}
\usepackage{setspace}

\newtheorem{theorem}{Theorem}[section]

\newtheorem{lemma}[theorem]{Lemma}

\newtheorem{remark}[theorem]{Remark}

\newtheorem{proposition}[theorem]{Proposition}

\newcommand{\wdin}{\upscr{d}{in}}
\newcommand{\wdout}{\upscr{d}{out}}
\newcommand{\nout}{\upscr{N}{out}}
\newcommand{\nin}{\upscr{N}{in}}
\newcommand{\amin}{\subscr{a}{min}}
\newcommand{\wdoutmax}{\upsubscr{d}{max}{out}}
\newcommand{\Nin}{\upscr{\mathcal{N}}{in}}

\newcommand{\Dout}{\subscr{\mathsf{D}}{out}}

\newcommand{\Adj}{\mathsf{A}}
\newcommand{\Lap}{\mathsf{L}}
\newcommand{\vertices}{\mathcal{V}}
\newcommand{\edges}{\mathcal{E}}

\newcommand{\Bgraph}{\mathcal{G}}

\newcommand{\GGa}{\GG}
\newcommand{\GGau}{\hat{\GG}}
\newcommand{\GGi}{\GGau_{\setminus\{4,11,25,45\}}} 
\newcommand{\GGf}{\GGau_{\setminus\{4,25,27\}}} 

\newcommand{\BB}{\mathcal{B}}
\renewcommand{\SS}{\mathcal{S}}

\newcommand{\UU}{\mathcal{U}}

\newcommand{\FF}{\mathcal{F}}
\newcommand{\GG}{\mathcal{G}}

\newcommand{\HH}{\mathcal{H}}

\newcommand{\OO}{\mathcal{O}}
\newcommand{\EE}{\mathcal{E}}
\newcommand{\MM}{\mathcal{M}}

\newcommand{\spn}{\mathrm{span}}

\newcommand{\ones}{\mathbf{1}}

\newcommand{\abs}[1]{\ensuremath{\left\lvert{#1}\right\rvert}}
\newcommand{\norm}[1]{\ensuremath{\| #1 \|}}

\newcommand{\Bnorm}[1]{\ensuremath{\Big \| #1 \Big \|}}

\newcommand{\real}{{\mathbb{R}}}
\newcommand{\realpositive}{{\mathbb{R}}_{>0}}
\newcommand{\realnonnegative}{{\mathbb{R}}_{\ge 0}}

\newcommand{\integerspositive}{\mathbb{Z}_{\geq 1}}

\newcommand{\eps}{\epsilon}

\newcommand{\until}[1]{\{1,\dots,#1\}}
\newcommand{\map}[3]{#1:#2 \rightarrow #3}
\newcommand{\setmap}[3]{#1:#2 \rightrightarrows #3}

\newcommand{\rrarrows}{\rightrightarrows}

\newcommand{\setdef}[2]{\{#1 \; | \; #2\}}

\newcommand{\Lie}{\mathcal{L}}
\newcommand{\SetLie}{{\mathcal{L}}}
\newcommand{\gradient}{\nabla}
\newcommand{\Eq}[1]{\mathrm{Eq}(#1)}

\newcommand{\ED}{{\rm ED}\xspace}

\newcommand{\FFED}{\FF_{\mathrm{ED}}}

\newcommand{\cLapgradnon}{X_{\texttt{lm+}\Lap\partial}}
\newcommand{\rLapgradnon}{X_{\texttt{dac+}\Lap\partial}}

\newcommand{\dac}{\texttt{dac}\xspace}
\newcommand{\lmplap}{\texttt{lm+}$\Lap\partial$\xspace}
\newcommand{\dacplap}{\texttt{dac+}$\Lap\partial$\xspace}

\newcommand{\HHp}{\HH_{P_l}}
\newcommand{\HHz}{\HH_0}
\newcommand{\nm}{n_{\text{max}}}
\newcommand{\FFa}{\FF_{\mathrm{aug}}^*}
\newcommand{\bFFa}{\overline{\FF}_{\mathrm{aug}}}

\newcommand{\vb}{\bar{v}}

\newcommand{\idf}{\mathrm{id}_{54}}

\newcommand{\tinv}{\textsc{trajectory invariance}\xspace} 
\DeclareMathAlphabet{\mathpzc}{OT1}{pzc}{m}{it}

\newcommand\subscr[2]{#1_{\textup{#2}}}
\newcommand\upscr[2]{#1^{\textup{#2}}}
\newcommand\upsubscr[3]{#1_{\textup{#2}}^{\textup{#3}}}

\newcommand{\oprocendsymbol}{\hbox{$\bullet$}}
\newcommand{\oprocend}{\relax\ifmmode\else\unskip\hfill\fi\oprocendsymbol}

\newcommand{\longthmtitle}[1]{\mbox{}\textup{\textsl{(#1):}}}

\renewcommand{\theenumi}{(\roman{enumi})}

\begin{document}

%

\begin{frontmatter}

  \title{Initialization-free distributed coordination for economic
    dispatch under varying loads and generator
    commitment\thanksref{footnoteinfo}}

  \thanks[footnoteinfo]{A preliminary version of this work appears at
    the 2014 Allerton Conference on Communication, Control, and
    Computing, Monticello, Illinois.}

  \author[ucsd]{Ashish Cherukuri}\ead{acheruku@ucsd.edu} \qquad
  \author[ucsd]{Jorge Cort\'{e}s}\ead{cortes@ucsd.edu}

  \address[ucsd]{Department of Mechanical and Aerospace Engineering,
    University of California, San Diego, CA, 92093, USA}

  \begin{abstract}
    This paper considers the economic dispatch problem for a network
    of power generating units communicating over a strongly connected,
    weight-balanced digraph.  The collective aim is to meet a power
    demand while respecting individual generator constraints and
    minimizing the total generation cost.  We design a distributed
    coordination algorithm consisting of two interconnected dynamical
    systems. One block uses dynamic average consensus to estimate the
    evolving mismatch in load satisfaction given the generation levels
    of the units. The other block adjusts the generation levels based
    on the optimization objective and the estimate of the load
    mismatch.  Our convergence analysis shows that the resulting
    strategy provably converges to the solution of the dispatch
    problem starting from any initial power allocation, and therefore
    does not require any specific procedure for initialization.  We
    also characterize the algorithm robustness properties against the
    addition and deletion of units (capturing scenarios with
    intermittent power generation) and its ability to track
    time-varying loads.  Our technical approach employs a novel
    refinement of the LaSalle Invariance Principle for differential
    inclusions, that we also establish and is of independent interest.
    Several simulations illustrate our results.
  \end{abstract}
\end{frontmatter}

\section{Introduction}\label{se:intro}

The advent of renewable energy sources and their integration into
electricity grids is making power generation and distribution an
increasingly decentralized problem.  The large-scale and highly
dynamic nature of the resulting grid optimization problems makes
traditional centralized, top-down approaches impractical because they
rely on the assumption of a fixed, limited number of generation units.
To solve these problems efficiently, there is a need to design
distributed algorithms that can handle dynamic loads, provide
plug-and-play capabilities, are robust against transmission and
generation failures, and adequately preserve the privacy of the
entities involved.  These considerations motivate us to consider here
the design of distributed algorithmic solutions to the economic
dispatch (\ED) problem, where a group of power generators aims to meet
a power demand while minimizing the total generation cost (the
summation of individual costs) and respecting the individual
generators' capacity constraints.  We are interested in the synthesis
of strategies that solve the \ED problem starting from any initial
power allocation, can handle time-varying loads, and are robust
against intermittent power generation caused by unit addition and
deletion.

\emph{Literature review:} The \ED problem has been traditionally
solved in a centralized manner, see e.g.~\citep{BHC-SR:90} and
references therein. Since distributed algorithmic solutions to grid
optimization problems are envisioned as part of the future smart
grid~\citep{HF:10}, this has motivated the appearance of a number of
distributed strategies for the \ED problem in the literature. While
there exists a broad variety in the assumptions made, a majority of
the works rely on the specific form of the solutions of the
optimization problem and propose consensus-based algorithms. Many
works consider convex, quadratic cost functions for the power
generators and perform consensus over their incremental costs under
undirected~\citep{ZZ-MC:12,SK-GH:12} or
directed~\citep{ADDG-STC-CNH:12,GB-AD-FLL-DN-BT:14} communication
topologies.  Some works consider general convex cost functions, like
we do here, but either do not consider capacity constraints on the
generators~\citep{RM-SD-BBC:12}, assume the initial power allocation
to meet the total load~\citep{AC-JC:13-tcns,AP-NQ-KMP:14}, or require
feedback on the power mismatch from the shift in frequency due to
primary droop control~\citep{WZ-WL-XW-LL-FF:14}.  Along with load and
capacity constraints,~\cite{GB-AD-FLL-DN-BT:14,VL-AV:13} consider
transmission losses, and~\cite{GB-AD-DN-BT-FLL:14} additionally take
into account valve-point loading effects and prohibited operating
zones. However, these constraints make the problem nonconvex and
prevent these works from obtaining theoretical guarantees on the
algorithm convergence properties.  In~\citep{LD-SG-RGH:12}, the
authors propose best-response dynamics for a potential-game
formulation of the nonconvex \ED problem, but the implementation
requires all-to-all communication among the
generators.~\cite{LX-SB:06,BJ-MJ:09} introduce distributed methods to
solve resource allocation problem very similar to the \ED problem, but
without taking into account individual agent constraints. Instead,
these are incorporated in the formulation of~\cite{AS-TK-MJ:12}, but
the proposed algorithm arrives at suboptimal solutions of the
optimization problem.
Our algorithm design and analysis rely on dynamic average consensus
and differential inclusions.  In dynamic average consensus, see
e.g.~\citep{RAF-PY-KML:06b,SSK-JC-SM:14-ijrnc} and references therein,
each agent has access to a time-varying input signal and interacts
with its neighbors in order to track the average of the input signals
across the network.  We build on our previous
work~\citep{AC-JC:13-tcns}, which requires a proper algorithm
initialization, and employ tools from dynamic average consensus to
synthesize a coordination strategy that converges from any initial
condition.
Regarding analysis, our technical approach builds on Lyapunov
stability tools for differential inclusions and nonsmooth systems, see
e.g.~\citep{AB-FC:99,JC:08-csm-yo} and references therein.  Of
particular importance is the work~\citep{AA-CE:10} for differential
equations, that provides a way to further refine the description of
omega-limit sets of trajectories by employing more than one
LaSalle-type function.

\emph{Statement of contributions:} We start with the formal definition
of the \ED problem for a network of power generators communicating
over a strongly connected, weight-balanced digraph.  The optimization
problem is convex as the individual cost functions are smooth and
convex, the load satisfaction is a linear constraint, and the
individual generators' capacities prescribe convex inequality
constraints. Our formulation is a simplification of the \ED problem in
its full generality, which in practice may have additional constraints
(e.g., transmission losses, line capacity constraints, valve-point
loading effects, ramp rate limits, prohibited operating zones) that
make it nonconvex. However, our developments show that obtaining a
provably correct algorithmic solution for the formulation here of the
\ED problem given our performance requirements (distributed,
convergent irrespective of initial condition, able to handle
time-varying loads, and robust to intermittent power generation) is
challenging.  Our first contribution is the design of a centralized
algorithm, termed ``load mismatch + Laplacian-nonsmooth-gradient''
dynamics, that solves the \ED problem starting from any initial power
allocation. This strategy has two components: one component seeks to
optimize the network generation cost while keeping constant the total
power generated; the other component is a feedback correction term
driven by the error between the desired total load and the network
generation. This latter term is responsible for ensuring that the
algorithm trajectories asymptotically satisfy the load satisfaction
constraint irrespective of the initial power allocation.  These
observations set the basis for our second contribution, which is the
synthesis of a distributed coordination algorithm, termed ``dynamic
average consensus + Laplacian-nonsmooth-gradient'' dynamics, with the
same convergence guarantees. Our design consists of two coupled
dynamical systems: a dynamic average consensus algorithm to estimate
the mismatch between generation and desired load in a distributed
fashion and a distributed Laplacian-nonsmooth-gradient dynamics that
employs these estimates to dynamically allocate the unit generation
levels.  The convergence analysis of both the centralized and
distributed algorithms relies on a combination of tools from algebraic
graph theory, nonsmooth analysis, set-valued dynamical systems, and
dynamic average consensus, and most notably on a refined version of
the LaSalle Invariance Principle for differential inclusions, which
constitutes our third contribution. Roughly speaking, the application
of the known LaSalle Invariance Principle would only establish
asymptotic convergence towards the network satisfaction of the total
load. Instead, the use of the refined version allows us, for each
algorithm, to establish global convergence of the trajectories to the
solutions of the \ED problem.  Our final contribution is the formal
characterization of the robustness properties of the distributed
algorithm. Building on the observation that the mismatch dynamics
between network generation and total load is exponentially convergent
and input-to-state stable, we establish the algorithm ability to track
time-varying loads and its robustness in scenarios with intermittent
power generation.



\section{Preliminaries}\label{sec:Prelim}

This section introduces basic concepts and preliminaries.
We begin with some notational
conventions. Let $\real$, $\realnonnegative$, $\realpositive$,
$\integerspositive$ denote the real, nonnegative real, positive real,
and positive integer numbers, resp. For $r \in \real$ we denote $\HH_r
= \setdef{x \in \real^n}{\ones_n^\top x = r}$.  The $2$- and
$\infty$-norms on $\real^n$ and their respective induced norms on
$\real^{n \times n}$ are denoted with $\norm{\cdot}$ and
$\norm{\cdot}_{\infty}$, resp. We let $B(x,\delta) = \setdef{y \in
  \real^n}{\norm{y-x} < \delta}$. For $D \subset \real^n$, $\overline
D$ denotes its closure.  For $x \in \real^n$, $x_i \in \real$ denotes
its $i$-th component. Given vectors $x, y \in \real^n$, $x \le y$ if
and only if $x_i \le y_i$ for all $i \in \until{n}$. We denote
$\ones_n = (1, \ldots, 1) \in \real^n$.A set-valued map
$\setmap{f}{\real^{n}}{\real^{m}}$ associates to each point in
$\real^{n}$ a set in $\real^{m}$. For a symmetric matrix $A \in
\real^{n \times n}$, $\lambda_{\min}(A)$ and $\lambda_{\max}(A)$
denote the minimum and maximum eigenvalues of $A$.
Finally, we let $[u]^{+} = \max \{0,u\}$ for $u \in \real$.

\subsection{Graph theory}\label{subse:Graph}

We present basic notions from algebraic graph theory
following~\citep{FB-JC-SM:08cor}.  A \emph{directed graph} (or
\emph{digraph}) is a pair $\Bgraph=(\vertices,\edges)$, with
$\vertices = \until{n}$ the vertex set and $ \edges \subseteq
\vertices\times \vertices $ the edge set.  A path is a sequence of
vertices connected by edges.  A digraph is \emph{strongly connected}
if there is a path between any pair of vertices. The sets of out- and
in-neighbors of $v$ are, resp., $\nout(v) = \setdef{w \in
  \vertices}{(v, w) \in \edges}$ and $\nin(v) = \setdef{w \in
  \vertices}{(w,v) \in \edges}$.  
A \emph{weighted digraph}
$\Bgraph=(\vertices,\edges,\Adj) $ is composed of a digraph $
(\vertices,\edges) $ and an \emph{adjacency matrix} $ \Adj \in
\mathbb{R}^{n\times n}_{\geq0} $ with $a_{ij}>0 $ iff $ (i,j)\in
\edges$. The weighted out- and in-degree of $i$ are, resp., $
\wdout(i) =\sum_{j=1}^{n}a_{ij} $ and $\wdin(i)=\sum_{j=1}^n a_{ji}$.
The \emph{Laplacian} matrix is $ \Lap = \Dout -\Adj$, where $ \Dout$
is the diagonal matrix with $(\Dout)_{ii}=\wdout(i) $, for all $ i \in
\{1,\ldots,n\}$.  Note that $\Lap\ones_n=0 $.  If $ \Bgraph $ is
strongly connected, then $0$ is a simple eigenvalue of~$\Lap$.
$\Bgraph$ is undirected if $ \Lap=\Lap^\top$. $\Bgraph$ is
\emph{weight-balanced} if $ \wdout(v) =\wdin(v) $, for all $ v \in
\vertices$ iff $ \ones_n^\top \Lap =0 $ iff $ \Lap+\Lap^\top \ge 0$.
Note that any undirected graph is weight-balanced.  If $\Bgraph $ is
weight-balanced and strongly connected, then $0$ is a simple
eigenvalue of $\Lap+\Lap^\top $. In such case, one has for $x \in
\real^n$,
\begin{equation}\label{eq:LapBound}
  x^\top (\Lap + \Lap^\top)x
  \ge   \lambda_2(\Lap + \Lap^\top) \Bnorm{x-\frac{1}{n}(\ones_n^\top
    x)\ones_n}^2  ,
\end{equation}
with $\lambda_2(\Lap + \Lap^\top\!)$ the smallest non-zero eigenvalue of
$\Lap + \Lap^\top\!$.

\subsection{Dynamic average consensus}\label{subse:dac}

Here, we introduce notions on dynamic average consensus
following~\citep{SSK-JC-SM:14-ijrnc}. Consider $n \in
\integerspositive$ agents communicating over a strongly connected,
weight-balanced digraph $\GG$ whose Laplacian is denoted as $\Lap$.
Each agent is associated with a state $x_i \in \real$ and an input
signal $t \mapsto u_i(t) \subset \real$ that is measurable and locally
essentially bounded. The aim is to provide a distributed dynamics such
that the state of each agent $x_i(t)$ tracks the average signal
$\frac{1}{n}\sum_{i = 1}^n u_i(t)$ asymptotically. This can be
achieved via the dynamics $\map{X_{\dac}}{\real^{2n}}{\real^{2n}}$,
\begin{subequations}
  \begin{align*}
    \dot x & = -\alpha x - \beta \Lap x - v + \nu u,
    \\
    \dot v & = \alpha \beta \Lap x,
  \end{align*}
\end{subequations}
where $\alpha , \beta, \nu > 0$ are design parameters and $v \in
\real^n$ is an auxiliary state.  If the initial condition satisfies
$\ones_n^\top v(0) = 0$ and the time-derivatives of the input signals
are bounded, then one can show, cf.~\citep[Corollary
4.1]{SSK-JC-SM:14-ijrnc}, that the error signal $t \mapsto \abs{x_i(t)
  -\frac{1}{n}\sum_{i = 1}^n u_i(t)}$ is ultimately bounded for each
$i \in \until{n}$. Moreover, this error vanishes if the input signal
converges to a constant value.

%

\subsection{Nonsmooth analysis and differential
  inclusions}\label{subse:nonsmooth}

We review here some notions from nonsmooth analysis and differential
inclusions following~\citep{JC:08-csm-yo}.  A function $f:\real^n
\rightarrow \real^m$ is \emph{locally Lipschitz} at $x \in \real^n$ if
there exist $L_x, \eps \in (0,\infty)$ such that $\norm{f(y) -
  f(y')} \le L_x\norm{y - y'}$, for all $y, y'\in B(x,\eps)$.  A
function $f:\real^n \rightarrow \real$ is \emph{regular} at $x \in
\real^n$ if, for all $v \in \real^n$, the right and generalized
directional derivatives of $f$ at $x$ in the direction of $v$
coincide, see~\citep{JC:08-csm-yo} for definitions of these notions.  A
function that is continuously differentiable at $x$ is regular at
$x$. Also, a convex function is regular.  A set-valued map
$\HH:\real^n \rrarrows \real^n$ is \emph{upper semicontinuous} at $x
\in \real^n$ if, for all $\eps \in (0,\infty)$, there exists $\delta
\in (0,\infty)$ such that $\HH(y) \subset \HH(x) + B(0,\eps)$ for all
$y \in B(x,\delta)$. Also, $\HH$ is \emph{locally bounded} at $x \in
\real^n$ if there exist $\eps, \delta \in (0,\infty)$ such that
$\norm{z} \le \eps$ for all $z \in \HH(y)$ and $y \in B(x,\delta)$.

Given a locally Lipschitz function $f:\real^n \rightarrow \real$, let
$\Omega_f$ be the set (of measure zero) of points where $f$ is not
differentiable. The \emph{generalized gradient} $\partial f: \real^n
\rrarrows \real^n$ is
\begin{equation*}
  \partial f(x) = \mathrm{co} \setdef{ \lim_{i \rightarrow \infty} 
    \gradient f(x_i)}{ x_i \rightarrow x, x_i \notin S \cup \Omega_f},
\end{equation*}
where $\mathrm{co}$ denotes convex hull and $S \subset \real^n$ is any
set of measure zero.  The map $\partial f$ is locally bounded, upper
semicontinuous, and takes non-empty, compact, and convex values. A
\emph{critical point} $x$ of $f$ satisfies~$0 \in \partial f(x)$.

Given a set-valued map $\HH: \real^n \rrarrows \real^n$, a
differential inclusion on $\real^n$~is
\begin{equation}\label{eq:ddsys}
  \dot x \in \HH(x) .
\end{equation}
A solution of~\eqref{eq:ddsys} on $[0,T] \subset \real$ is an
absolutely continuous map $x:[0,T]\rightarrow \real^n$ that
satisfies~\eqref{eq:ddsys} for almost all $t \in [0,T]$.  If $\HH$ is
locally bounded, upper semicontinuous, and takes non-empty, compact,
and convex values, then existence of solutions is guaranteed. The set
of equilibria of~\eqref{eq:ddsys} is $\Eq{\HH} = \setdef{x \in
  \real^n}{0 \in \HH(x) }$.  A set $ S \subset \real^n $ is
\emph{weakly (resp., strongly) positively invariant}
under~\eqref{eq:ddsys} if, for each $x \in S$, at least a solution
(resp., all solutions) starting from~$x$ is (resp., are) entirely
contained in $S$.  For dynamics with uniqueness of solution, both
notions coincide and are referred as \emph{positively invariant}.
Given a locally Lipschitz function $f: \real^n \rightarrow \real$, the
\emph{set-valued Lie derivative} $\SetLie_{\HH}f: \real^n \rrarrows
\real$ of $f$ with respect to~\eqref{eq:ddsys}  is
\begin{align*}
  \SetLie_{\HH}f (x) = \setdef{a \in \real}{\exists v \in \HH(x)
    \text{ s.t. } \zeta^\top v=a & \text{ for all }
    \\
    & \zeta \in \partial f(x)}.
\end{align*}
For a trajectory $t \mapsto \varphi(t)$, $\varphi(0) \in \real^n$
of~\eqref{eq:ddsys}, the evolution of $f$ along it satisfies 
\begin{align*}
  \frac{d}{dt} f(\varphi(t)) \in \SetLie_{\HH} f (\varphi(t))
\end{align*}
for almost all $t \ge 0$. The omega-limit set of the trajectory,
denoted $\Omega(\varphi)$, is the set of all points $y \in \real^n$
for which there exists a sequence $\{t_k\}_{k=1}^{\infty}$ with $t_k
\to \infty$ and $\lim_{k \to \infty} \varphi(t_k) = y$.  If the
trajectory is bounded, then the omega-limit set is nonempty, compact,
connected, and weakly invariant. These tools allow us to characterize
the asymptotic behavior of solutions of differential inclusions. In
Appendix~\ref{app:refined-LaSalle} we develop a novel refinement of
the LaSalle Invariance Principle for differential inclusions, see
e.g.,~\citep{JC:08-csm-yo}, which is suitable for the analysis of the
coordination algorithms.

\section{Problem statement}\label{sec:problem}

This section presents the network model and the economic dispatch
problem we set out to solve in a distributed and robust fashion.
Consider $n \in \integerspositive$ power generators communicating over
a strongly connected and weight-balanced digraph
$\Bgraph=(\vertices,\edges,\Adj)$. Each generator corresponds to a
vertex in the digraph and an edge $(i,j)$ represents the ability of
generator~$j$ to send information to generator $i$. The cost of power
generation for unit $i$ is measured by $f_i:\real \rightarrow
\realnonnegative$, assumed to be convex and continuously
differentiable.  Representing the power generated by unit $i$ by $P_i
\in \real$, the total cost incurred by the network with the power
allocation $P = (P_1, \dots, P_n) \in \real^n$ is measured by
$f:\real^n \rightarrow \realnonnegative$ as
\begin{align*}
  f(P) = \sum_{i=1}^{n}f_i(P_i) .
\end{align*}
Note that $f$ is convex and continuously differentiable.  The
generators aim to minimize the total cost $f(P)$ while meeting the
total power load $P_l \in \realpositive$, i.e., $\sum_{i=1}^n P_i =
P_l$.  Each generator has an upper and a lower limit on the power it
can produce, $P_i^m \le P_i \le P_i^M$ for $i \in \until{n}$.
Formally, the economic dispatch (\ED) problem is
\begin{subequations}\label{eq:EDP}
  \begin{align}
    \mathrm{minimize} \quad & f(P), \label{eq:conobjective}
    \\
    \text{subject to} \quad & \ones_n^\top P =
    P_l, \label{eq:equalitycons}
    \\
    & P^{m} \le P \le P^{M}. \label{eq:inequalitycons}
  \end{align}
\end{subequations}
The constraint~\eqref{eq:equalitycons} is the \emph{load condition}
and~\eqref{eq:inequalitycons} are the \emph{box constraints}. The set
of allocations satisfying the box constraints is $\FF_B = \setdef{P
  \in \real^n}{P^m \le P \le P^M}$. Further, we denote the feasibility
set of~\eqref{eq:EDP} as $\FFED = \FF_B \cap \HHp = \setdef{P \in
  \real^n}{P^m \le P \le P^M \text{ and } \ones_n^\top P = P_l}$ and
the set of solutions as $\FFED^{*}$. Since $\FFED$ is compact,
$\FFED^*$ is compact.
Note that $P^M \in \FFED$ implies $\FFED = \{P^M\}$. Similarly $P^m
\in \FFED$ implies $\FFED = \{P^m\}$.  Therefore, we assume $P^M$ and
$P^m$ are not feasible.

Our objective is to design a distributed coordination algorithm that
allows the team of generators to solve the \ED problem~\eqref{eq:EDP}
starting from any initial condition, can handle time-varying loads,
and is robust to intermittent power generation.  

\begin{remark}\longthmtitle{Additional practical constraints}
  \label{re:load-agg-1}
  {\rm We do not consider here, for simplicity, other constraints on
    the \ED problem such as transmission losses, transmission line
    capacities, valve-point loading effects, ramp rate limits, and
    prohibited operating zones.  As our forthcoming treatment will
    show, the design and analysis of algorithmic solutions to the \ED
    problem without these additional constraints is already quite
    challenging given our performance requirements.
    Nevertheless, Remark~\ref{re:load-agg-2} later comments on how to
    adapt our algorithm to deal with more general scenarios.
    \oprocend }
\end{remark}

Our design strategy relies on the following reformulation of the \ED
problem without inequality constraints.  Consider the modified \ED
problem
\begin{subequations}\label{eq:modified-edp}
  \begin{align}
    \mathrm{minimize} \quad & f^{\eps}(P), 
    \\
    \text{subject to} \quad & \ones_n^\top P = P_l,
  \end{align}
\end{subequations}
where the objective function is
\begin{align*}
  f^{\eps}(P) = \sum_{i=1}^n f_i(P_i) +
  \frac{1}{\eps}(\sum_{i=1}^n([P_i-P_i^M]^+ + [P_i^m - P_i]^+)).
\end{align*}
This corresponds to each generator $i \in \until{n}$ having the
modified local cost
\begin{align*}
  f_i^{\eps} (P_i) = f_i(P_i) + \frac{1}{\eps}([P_i - P_i^M]^+ +
  [P_i^m -  P_i]^+).
\end{align*}
Note that $f_i^{\eps}$ is convex, locally Lipschitz, and continuously
differentiable on $\real$ except at $P_i = P_i^m$ and $P_i = P_i^M$.
%
%
Moreover, the total cost $f^{\eps}$ is convex, locally Lipschitz, and
regular.
According to our previous work~\citep[Proposition~5.2]{AC-JC:13-tcns},
the solutions to the original~\eqref{eq:EDP} and the
modified~\eqref{eq:modified-edp} \ED problems coincide for $\eps \in
\realpositive$ such that
\begin{equation}\label{eq:EpsBound}
  \eps < \frac{1}{2\max_{P \in \FFED}\norm{\gradient
      f(P)}_{\infty}}.
\end{equation}
Throughout the paper, we assume the parameter $\eps$ satisfies this
condition.  A useful fact is that $P^* \in \real^n$ is a solution
of~\eqref{eq:modified-edp} if and only if there exists $\mu \in \real$
such that
\begin{align}\label{eq:modified-edp-sol}
  \mu \ones_n \in \partial f^{\eps} (P^*) \quad \text{ and }
  \quad \ones_n^\top P^* = P_l.
\end{align}

\section{Robust centralized algorithmic
  solution}\label{sec:centralized}

This section presents a robust strategy to make the network power
allocation converge to the solution set of the \ED problem starting
from any initial condition.  Even though this algorithm is
centralized, its design provides enough insight to tackle later the
design of a distributed algorithmic solution.  Consider the ``load
mismatch + Laplacian-nonsmooth-gradient'' (abbreviated \lmplap)
dynamics, represented by the set-valued map $\setmap{\cLapgradnon}
{\real^n}{\real^n}$,
\begin{equation}\label{eq:central-lap-gradient} 
  \dot{P} \in  -\Lap \partial f^{\eps}(P) + \frac{1}{n}
  (P_l - \ones_n^\top P) \ones_n ,
\end{equation}
where $\Lap$ is the Laplacian associated to the strongly connected and
weight-balanced communication digraph~$\Bgraph$.  For each generator,
the first term seeks to minimize the total cost while leaving
unchanged the total generated power. The second term is a feedback
element that seeks to drive the units towards the satisfaction of the
load. The first term is computable using information from its
neighbors but the second term requires them to know the aggregated
state of the whole network, which makes it not directly implementable
in a distributed manner.  The next result shows that the trajectories
of~\eqref{eq:central-lap-gradient} converge to the set of solutions of
the \ED problem.

\begin{theorem}\longthmtitle{Convergence of the trajectories of
    $\cLapgradnon$ to the solutions of \ED
    problem}\label{th:central-convergence}
  The trajectories of~\eqref{eq:central-lap-gradient} starting
  from any point in $\real^n$ converge to the set of solutions
  of~\eqref{eq:EDP}.
\end{theorem}
\begin{pf}
  Our proof strategy proceeds by applying the refined LaSalle
  Invariance Principle for differential inclusions established in
  Appendix~\ref{app:refined-LaSalle},
  cf. Proposition~\ref{pr:refined-lasalle-nonsmooth}.  Consider the
  following function $\map{V_1}{\real^n}{\realnonnegative}$,
  \begin{equation*}
    V_1(P) = \frac{1}{2}(P_l - \ones_n^\top P)^2.
  \end{equation*}
  The set-valued Lie derivative of $V_1$ along $\cLapgradnon$ is
  \begin{align*}
    \SetLie_{\cLapgradnon}V_1(P) = \{-(P_l - \ones_n^\top P)^2\} = \{- 2
    V_1(P)\}.
  \end{align*}
  Thus, starting at any $P(0) \in \real^n$, the trajectory of
  $\cLapgradnon$ satisfies $V_1(P(t)) = V_1(P(0)) e^{-2t}$ and its
  omega-limit set (provided the trajectory is bounded, a fact that we
  assume is true for now and establish later) is contained in~$\HHp$.
  In the notation of Proposition~\ref{pr:refined-lasalle-nonsmooth},
  $\HH_{P_l}$ plays the role of the closed submanifold $\SS$
  of~$\real^n$.  We next show that the hypotheses of this result hold.
  In the notation of the Lemma~\ref{le:continuity-set-valued-gen}, the
  function $f^{\eps}$, the map $(P,\zeta) \mapsto - \Lap \zeta$, and
  the set-valued map $P \rrarrows -\Lap \partial f^{\eps}(P)$ play the
  role of $W$, $g$, and $F$, respectively (our choice of $F$ is
  because the dynamics $\cLapgradnon$ takes the form $\dot P \in -
  \Lap \partial f^{\eps}(P)$ on~$\SS=\HH_{P_l}$).  Notice that $\zeta
  \mapsto -\Lap \zeta$ is a continuous map and, since $\GG$ is
  strongly connected and weight-balanced, we have $\zeta^\top (-\Lap
  \zeta) = -\frac{1}{2}\zeta^\top (\Lap + \Lap^\top) \zeta \le 0$ for
  any $\zeta \in \partial f^{\eps}(P)$. Therefore,
  Lemma~\ref{le:continuity-set-valued-gen}\ref{as:continuity-1} is
  satisfied.  Moreover, if $\zeta^\top (-\Lap \zeta) = 0$ for some
  $\zeta \in \partial f^{\eps} (P)$, then $\zeta \in
  \spn\{\ones_n\}$. Since for $P \in \HH_{P_l}$, we have
  \begin{equation*}
    \SetLie_{\cLapgradnon} f^{\eps}(P) = \setdef{-\zeta^\top \Lap
      \zeta}{\zeta \in \partial f^{\eps}(P)}, 
  \end{equation*}
  we deduce $0 \in \SetLie_{\cLapgradnon} f^{\eps}(P)$, i.e.,
  Lemma~\ref{le:continuity-set-valued-gen}\ref{as:continuity-2} holds.
  The application of Lemma~\ref{le:continuity-set-valued-gen} then
  yields that
  Proposition~\ref{pr:refined-lasalle-nonsmooth}\ref{as:refined-lasalle-2}
  holds too.  In addition, from the above analysis, note that if $0
  \in \SetLie_{\cLapgradnon} f^{\eps}(P)$ for some $P \in \HH_{P_l}$,
  then there exists $\mu \in \real$ such that $\mu \ones_n
  \in \partial f^{\eps}(P)$ and, from~\eqref{eq:modified-edp-sol}, $P$
  is a solution of~\eqref{eq:EDP}.  Therefore, $\setdef{P \in \HHp}{0
    \in \SetLie_{\cLapgradnon}f^{\eps} (P)}$ is the set of solutions
  of the \ED problem and belongs to a level set of~$f^{\eps}$, which
  establishes that
  Proposition~\ref{pr:refined-lasalle-nonsmooth}\ref{as:refined-lasalle-1}
  also holds.

  To be able to apply Proposition~\ref{pr:refined-lasalle-nonsmooth}
  and conclude the proof, it remains to show that the trajectories
  of~$\cLapgradnon$ are bounded. We reason by contradiction, i.e.,
  assume there exists a trajectory $t \mapsto P(t)$, $P(0) \in
  \real^n$ of $\cLapgradnon$ such that $\norm{P(t)} \to \infty$. From
  the analysis above, we know that along this trajectory $\ones_n^\top
  P(t) \to P_l$ and $f^{\eps}(P(t)) \to \infty$ (as $f^{\eps}$ is
  radially unbounded). Therefore, there exist a sequence of times
  $\{t_k\}_{k=1}^{\infty}$ with $t_k \to \infty$ such that for all $k
  \in \integerspositive$,
  \begin{equation}\label{eq:unbounded-cond1}
    \abs{\ones_n^\top P(t_k) - P_l}  \!<\! \frac{1}{k}  \; \text{and} \; 
    \max \SetLie_{\cLapgradnon} f^{\eps} (P(t_k)) \!>\! 0.
  \end{equation}
  This implies that there exists a sequence $\{\zeta_k\}_{k =
    1}^{\infty}$ with $\zeta_k \in \partial f^{\eps} (P(t_k))$ such
  that, for all $k \in \integerspositive$,
  \begin{align}
    &- \zeta_k^\top \Lap \zeta_k + \frac{1}{n}(P_l - \ones_n^\top
    P(t_k)) (\ones_n^\top \zeta_k) > 0 \notag
    \\
    &\Rightarrow -\zeta_k^\top \Bigl( \frac{\Lap + \Lap^\top}{2} \Bigr)
    \zeta_k + \frac{1}{nk} \abs{\ones_n^\top \zeta_k} >
    0 \label{eq:unbounded-cond2}
    \\
    &\Rightarrow -\frac{\lambda_2(\Lap + \Lap^\top)}{2} \Bnorm{\zeta_k -
      \frac{1}{n} (\ones_n^\top \zeta_k) \ones_n}^2 + \frac{1}{nk}
      \abs{\ones_n^\top \zeta_k} > 0, \notag
  \end{align}
  where we have used~\eqref{eq:unbounded-cond1} in the first
  implication and~\eqref{eq:LapBound} in the second.  Next, we
  consider two cases depending on whether (a) $\abs{\ones_n^\top
    \zeta_k}$ is bounded or (b) $\abs{\ones_n^\top \zeta_k} \to
  \infty$. In case (a), taking the limit $k \to \infty$ in the last
  inequality of~\eqref{eq:unbounded-cond2}, we get
  \begin{equation}\label{eq:limit-cond}
    \lim_{k \to \infty} \Bnorm{\zeta_k -
      \frac{1}{n} (\ones_n^\top \zeta_k) \ones_n} = 0.
  \end{equation}
  Since, $\norm{P(t)} \to \infty$ and $\ones_n^\top P(t) \to P_l$,
  there exist $i,j \in \until{n}$ such that $P_i(t_k) \to \infty$ and
  $P_j(t_k) \to -\infty$. Let $P^* \in \FFED^*$ and $\mu \ones_n \in
  \partial f^{\eps} (P^*)$ for some $\mu \in \real$. Then, without
  loss of generality, we assume that $P^*_i \le P_i(t_k) \le
  P_i(t_{k+1})$ and $P^*_j \ge P_j(t_k) \ge P_j(t_{k+1})$ for all $k$.
  This fact along with the expression of $\setmap{\partial
    f_i^{\eps}}{\real}{\real}$,
 \begin{align*}
   \partial f_i^{\eps}(P_i) =
   \begin{cases}
     \{ \gradient f_i(P_i) - \frac{1}{\eps} \} & \quad \text{if }
     P_i <P_i^m,
     \\
     [\gradient f_i(P_i) - \frac{1}{\eps}, \gradient f_i(P_i)] &
     \quad \text{if } P_i=P_i^m ,
     \\
     \{ \gradient f_i(P_i) \} & \quad \text{if } P_i^m<P_i<P_i^M ,
     \\
     [\gradient f_i(P_i), \gradient f_i(P_i)+ \frac{1}{\eps}] &
     \quad \text{if } P_i=P_i^M ,
     \\
     \{ \gradient f_i(P_i)+ \frac{1}{\eps} \} & \quad \text{if }
     P_i > P_i^M.
   \end{cases}
  \end{align*}
  gives us the following property for all $k \in \integerspositive$,
  \begin{subequations}
    \begin{align}
      \min \partial f_i^{\eps} (P_i(t_k)) \ge \mu,  & \quad 
      \max \partial f_j^{\eps} (P_j(t_k)) \le \mu,
      \label{eq:unbounded-cond3-a}
      \\
      \min \partial f_i^{\eps}(P_i(t_{k+1})) & \ge \max \partial
      f_i^{\eps} (P_i(t_k)), \label{eq:unbounded-cond3-b}
      \\
      \max \partial f_j^{\eps}(P_j(t_{k+1})) & \le \min \partial
      f_j^{\eps} (P_j(t_k)). \label{eq:unbounded-cond3-c}
    \end{align}
  \end{subequations}
  Note that the limit~\eqref{eq:limit-cond} yields $\lim_{k \to
    \infty} \abs{(\zeta_k)_i - (\zeta_k)_j} = 0$.  On the other hand,
  from~\eqref{eq:unbounded-cond3-b}-\eqref{eq:unbounded-cond3-c}, we
  obtain $\abs{(\zeta_k)_i - (\zeta_k)_j} \le \abs{(\zeta_{k+1})_i -
    (\zeta_{k+1})_j}$ for all $k$.  Therefore, we obtain $(\zeta_k)_i
  = (\zeta_k)_j$ for all $k$ and from~\eqref{eq:unbounded-cond3-a}, we
  get $\mu = (\zeta_k)_i = (\zeta_k)_j$ for all $k$.
  From~\eqref{eq:unbounded-cond3-b}-\eqref{eq:unbounded-cond3-c}, this
  further implies that $\mu \in \partial f_i^{\eps} (x)$ for all $x
  \in [P^*_i , \infty)$ and that $\mu \in \partial f_j^{\eps} (x)$ for
  all $x \in (-\infty,P^*_j]$. Using this fact, one can construct an
  unbounded set of solutions to the \ED problem in the following
  manner. First, fix all the components of $P^*$ except $i$ and $j$.
  Now pick any $x \in \realnonnegative$ and consider $P_i^* + x$ and
  $P_j^* -x$. From what we have reasoned so far, all such points that
  we obtain by varying $x$ are solutions to the \ED problem as they
  satisfy~\eqref{eq:modified-edp-sol}.  This contradicts the fact that
  $\FFED^*$ is bounded.
  
  In case (b), assume without loss of generality that $\ones_n^\top
  \zeta_k \to \infty$ (the argument for $\ones_n^\top \zeta_k \to
  -\infty$ follows similarly).  As reasoned above, there exists $j \in
  \until{n}$ such that $P_j(t_k) \to - \infty$ and there exists $\mu
  \in \real$ such that $(\zeta_k)_j \le \mu$ for all $k \in
  \integerspositive$. Using this fact, we upper bound the left hand
  side of the inequality~\eqref{eq:unbounded-cond2} by
  \begin{align}
    & -\frac{\lambda_2(\Lap + \Lap^\top)}{2}  \Bnorm{\zeta_k -
      \frac{1}{n} (\ones_n^\top \zeta_k) \ones_n}^2 + 
      \frac{1}{nk} (\ones_n^\top \zeta_k) \notag
      \\
    & \le -\frac{\lambda_2(\Lap + \Lap^\top)}{2} \Bigl( (\zeta_k)_j - 
    \frac{1}{n} (\ones_n^\top \zeta_k) \Bigr)^2 + 
      \frac{1}{nk} (\ones_n^\top \zeta_k) \notag
      \\
    & \le -\frac{\lambda_2(\Lap + \Lap^\top)}{2} \Bigl( \mu - 
    \frac{1}{n} (\ones_n^\top \zeta_k) \Bigr)^2 + 
    \frac{1}{nk} (\ones_n^\top \zeta_k), \label{eq:unbounded-cond4}
  \end{align}
  where the last inequality is valid for all but a finite number of
  $k$.  Hence, as $\ones_n^\top \zeta_k \to \infty$, there is $\bar{k}
  \in \integerspositive$ such that the expression
  in~\eqref{eq:unbounded-cond4} is negative for~$k \ge \bar{k}$,
  contradicting~\eqref{eq:unbounded-cond2}. Thus, we conclude the
  trajectories are~bounded. \qed
\end{pf}

From the proof above, it is interesting to note that the feedback
term~\eqref{eq:central-lap-gradient} drives the mismatch between
generation and load to zero at an exponential rate, no matter what the
initial power allocation. This is a good indication of its robustness
properties: time-varying loads or scenarios with generators going down
and coming back online can be handled as long as the rate of these
changes is lower than the exponential rate of convergence associated
to the load satisfaction.
We provide a formal characterization of these properties for the
distributed implementation of this strategy in the next section.

\section{Robust distributed algorithmic
  solution}\label{sec:robust-dist-edp}

This section presents a distributed strategy to solve the \ED problem
starting from any initial power allocation. We build on the
centralized design presented in Section~\ref{sec:centralized}. We also
formally characterize the robustness properties against addition
and deletion of generators and time-varying loads.

Given the discussion on the centralized nature of the
dynamics~\eqref{eq:central-lap-gradient}, the core idea of our design
is to employ a dynamic average consensus algorithm that allows each
unit in the network to estimate the mismatch in load satisfaction.  To
this end, we assume the total load $P_l$ is only known to one
generator $r \in \until{n}$ (its specific identity is
arbitrary). Following Section~\ref{subse:dac}, consider the dynamics,
\begin{subequations}
  \begin{align*}
    \dot z & = -\alpha z - \beta \Lap z - v + \nu_2 (P_l e_r -P),
    \\
    \dot v & = \alpha \beta \Lap z, 
  \end{align*}
\end{subequations}
where $e_r \in \real^n$ is the unit vector along the $r$-th direction
and $\alpha, \beta, \nu_2 >0$ are design parameters.
Note that this dynamics is distributed over the communication graph
$\Bgraph$. For each $i \in \until{n}$, $z_i$ plays the role of an
estimator associated to $i$ which aims to track the average signal $t
\mapsto \frac{1}{n}(P_l - \ones_n^\top P(t))$. This observation
justifies substituting the feedback term
in~\eqref{eq:central-lap-gradient} by $z \in \real^n$, giving rise to
the ``dynamic average consensus + Laplacian-nonsmooth-gradient''
dynamics, abbreviated \dacplap for convenience, mathematically
represented by the set-valued map
$\setmap{\rLapgradnon}{\real^{3n}}{\real^{3n}}$,
\begin{subequations}\label{eq:robust-edp-dist}
  \begin{align}
    \dot P & \in -\Lap \partial f^{\eps}(P) + \nu_1
    z, \label{eq:red-one}
    \\
    \dot z & = -\alpha z - \beta \Lap z - v + \nu_2 (P_l e_r -P),
    \label{eq:red-two}
    \\
    \dot v & = \alpha \beta \Lap z, \label{eq:red-three}
  \end{align}
\end{subequations}
where $\nu_1 > 0$ is a design parameter.
Unlike~\eqref{eq:central-lap-gradient}, this dynamics is distributed,
as each agent only needs to interact with its neighbors to implement
it.

\subsection{Convergence analysis}\label{sec:convergence-distributed}

Here we characterize the asymptotic convergence properties of the
\dacplap dynamics. We start by establishing
an important fact on the omega-limit set of any trajectory
of~\eqref{eq:robust-edp-dist} with initial condition in $ \real^n
\times \real^n \times \HH_0$.

\begin{lemma}\longthmtitle{Characterizing the omega-limit set of the
    trajectories of the \dacplap
    dynamics}\label{le:omega-limit-set-edp}
  The omega-limit set of any trajectory of~\eqref{eq:robust-edp-dist}
  with initial condition $(P_0,z_0,v_0) \in \real^n \times \real^n
  \times \HH_0$ is contained in $\HHp \times \HHz \times \HHz$.
\end{lemma}
\begin{pf}
  From~\eqref{eq:red-three}, note that $\ones_n^\top \dot v =
  0$. Since $v_0 \in \HH_0$, this implies that $\ones_n^\top v(t) = 0$
  for all $t \ge 0$. Now, define $\zeta(t) = \ones_n^\top P(t) - P_l$
  and note that
  \begin{align*}
    \dot \zeta(t) = \ones_n^\top \dot P(t) = \nu_1 \ones_n^\top z(t) ,
  \end{align*}
  where we have used~\eqref{eq:red-one}, and
  \begin{align*}
    \ddot \zeta(t) & = \nu_1 \ones_n^\top \dot z(t)
    \\
    & = \nu_1 \ones_n^\top (-\alpha z(t) - \beta \Lap z(t) - v(t) +
    \nu_2 (P_l e_k - P(t))
    \\
    & = - \nu_1 \alpha (\ones_n^\top z(t)) - \nu_1 \nu_2 \zeta(t) = -
    \alpha \dot \zeta(t) - \nu_1 \nu_2 \zeta(t) ,
  \end{align*}
  where we have used~\eqref{eq:red-two}.  We write this system as a
  first-order one by defining $x_1 = \zeta$ and $x_2 = \dot \zeta$ to 
  get
  \begin{equation}\label{eq:mismatch-dyn}
    \begin{bmatrix}
      \dot x_1
      \\
      \dot x_2
    \end{bmatrix}
    = 
    \begin{bmatrix}
      0 & 1
      \\
      -\nu_1 \nu_2  & -\alpha
    \end{bmatrix}
    \begin{bmatrix}
      x_1
      \\
      x_2
    \end{bmatrix}.
  \end{equation}
  Evaluating the Lie derivative of the positive definite, radially
  unbounded function $V_2(x_1,x_2) = \nu_1 \nu_2 x_1^2 + x_2^2$ along
  the above dynamics and applying the LaSalle Invariance
  Principle~\citep{HKK:02}, we deduce that $x_1(t) \to 0$ and $x_2(t)
  \to 0$ as $t \to \infty$, that is, $\ones_n^\top P(t) \to P_l$ and
  $\ones_n^\top z(t) \to 0$. Since the system~\eqref{eq:mismatch-dyn}
  is linear, the convergence is exponential. \qed
\end{pf}

The next result builds on this fact and
Proposition~\ref{pr:refined-lasalle-nonsmooth} to establish that the
trajectory of power allocations under~\eqref{eq:robust-edp-dist}
converges to the solution set of the $\ED$ problem.

\begin{theorem}\longthmtitle{Convergence of the
    \dacplap dynamics to the solutions of
    \ED problem}\label{th:convergence-rlngd}
  For $\alpha, \beta, \nu_1, \nu_2>0$ with
  \begin{equation}\label{eq:alpha-beta-cond-n}
    \frac{\nu_1}{\beta \nu_2 \lambda_2(\Lap + \Lap^\top)} +
    \frac{\nu_2^2 \lambda_{\max}(\Lap^\top \Lap)}{2 \alpha} <
    \lambda_2(\Lap + \Lap^\top) ,
  \end{equation}
  the trajectories of~\eqref{eq:robust-edp-dist} starting from any
  point in $\real^n \times \real^n \times \HH_0$ converge to the set
  $\FFa = \setdef{(P,z,v) \in \FFED^* \times \{0\} \times \real^n}{v=
    \nu_2 (P_l e_r -P)}$.
\end{theorem}
\begin{pf}
  Our proof strategy is based on the refined LaSalle Invariance
  Principle for differential inclusions established in
  Appendix~\ref{app:refined-LaSalle},
  cf. Proposition~\ref{pr:refined-lasalle-nonsmooth}. Before
  justifying that all its hypotheses are satisfied, we reformulate the
  expression for the dynamics to help simplify the analysis.
  Consider first the change of coordinates, $(P,z,v) \mapsto
  (P,z,\vb)$, with $\vb = v - \nu_2 (P_l e_r - P)$.  The set-valued
  map $\rLapgradnon$ then takes the form
  \begin{align*}
    &\rLapgradnon (P,z,\vb) = \setdef{( - \Lap \zeta + \nu_1 z,
      -(\alpha + \beta \Lap) z -\vb,
      \\
      & \qquad \qquad (\alpha \beta \Lap + \nu_1 \nu_2) z - \nu_2 \Lap
      \zeta) \in \real^{3n}}{ \zeta \in \partial f^{\eps}(P)} .
  \end{align*}
  The change of coordinates shifts the equilibrium of the consensus
  dynamics to the origin.  Under the additional change of coordinates
  $(P,z,\vb) \mapsto (P,\xi_1,\xi_2)$, with
  \begin{equation}\label{eq:zv2xi}
    \begin{bmatrix}
      \xi_1
      \\ 
      \xi_2
    \end{bmatrix}
    = 
    \begin{bmatrix}
      I & 0
      \\ 
      \alpha I & I
    \end{bmatrix}
    \begin{bmatrix}
      z
      \\
      \vb
    \end{bmatrix}, 
  \end{equation}
  the set-valued map $\rLapgradnon$ takes the form
  \begin{align}
    \rLapgradnon &(P,\xi_1,\xi_2) = \setdef{( - \Lap \zeta + \nu_1
      \xi_1, -\beta \Lap \xi_1 -\xi_2, \label{eq:rewrite-edp2-2}
      \\
      & \nu_1 \nu_2 \xi_1 -\alpha \xi_2 - \nu_2 \Lap \zeta) \in
      \real^{3n}}{ \zeta \in \partial f^{\eps}(P)} . \notag
  \end{align}
  This extra change of coordinates makes it easier to~identify the
  candidate Lyapunov
  function~$\map{V_3}{\real^{3n}}{\realnonnegative}$,
  \begin{align*}
    V_3(P,\xi_1,\xi_2) = f^{\eps}(P) + \frac{1}{2}(\nu_1 \nu_2
    \norm{\xi_1}^2 + \norm{\xi_2}^2) .
  \end{align*}
  For convenience, denote the overall change of coordinates by
  $\map{D}{\real^{3n}}{\real^{3n}}$,
  \begin{align*}
    (P,\xi_1,\xi_2) = D(P,z,v) = (P,z,v +\alpha z- \nu_2 (P_l e_r -P)) .
  \end{align*}
  Our analysis now focuses on proving that, in the new coordinates,
  the trajectories of~\eqref{eq:robust-edp-dist} converge to the set
  \begin{align*}
    \bFFa & = D(\FFa) = \FFED^* \times \{0\} \times \{0\}.
  \end{align*}
  Note that $D(\HH_{P_l} \times \HH_0 \times \HH_0) = \HH_{P_l} \times
  \HH_0 \times \HH_0$ and therefore, from
  Lemma~\ref{le:omega-limit-set-edp}, the omega-limit set of a
  trajectory $t \mapsto (P(t),\xi_1(t),\xi_2(t))$ starting in
  $D(\real^{n} \times \real^{n} \times \HH_0)$ belongs to~$\HH_{P_l}
  \times \HH_0 \times \HH_0$.

  Our next step is to show that the hypotheses of
  Proposition~\ref{pr:refined-lasalle-nonsmooth} are satisfied where
  $\HH_{P_l} \times \HH_0 \times \HH_0$ and $V_3$ play the role of the
  closed submanifold $\SS$ of~$\real^{3n}$ and the function~$W$,
  respectively.  To do so, we resort to
  Lemma~\ref{le:continuity-set-valued-gen}.  Define the continuous
  function $\map{g}{\real^{3n} \times \real^{3n}}{\real^{3n}}$ by
  \begin{multline*}
    g (P,\xi_1,\xi_2,\hat{\zeta}) = ( - \Lap \hat{\zeta}_1 + \nu_1
    \xi_1, -\beta \Lap \xi_1 -\xi_2,
    \\
    \nu_1 \nu_2 \xi_1 -\alpha \xi_2 -
    \nu_2 \Lap \hat{\zeta}_1) ,
  \end{multline*}
  and note that the dynamics~\eqref{eq:rewrite-edp2-2} can be
  expressed as $\rLapgradnon(P,\xi_1,\xi_2) =
  \setdef{g(P,\xi_1,\xi_2,\hat{\zeta})}{\hat{\zeta} \in \partial
    V_3(P,\xi_1,\xi_2)}$.  For $(P,\xi_1,\xi_2) \in \HH_{P_l} \times
  \HH_0 \times \HH_0$ and $\hat{\zeta} \in \partial
  V_3(P,\xi_1,\xi_2)$,
  \begin{align}\label{eq:liederivative-element}
    \hat{\zeta}^\top g(P,\xi_1,\xi_2,\hat{\zeta}) & = - \zeta^\top
    \Lap \zeta + \nu_1 \zeta^\top \xi_1 - \beta \nu_1 \nu_2 \xi_1^\top
    L \xi_1 \notag
    \\
    & \qquad - \alpha \norm{\xi_2}^2 - \nu_2 \xi_2^\top L \zeta,
  \end{align}
  where we have used that $\zeta = \hat{\zeta}_1 \in \partial
  f^{\eps}(P)$, $\hat{\zeta}_2 = \nu_1 \nu_2 \xi_1$, and
  $\hat{\zeta}_3 = \xi_2$.  Since the digraph $\GG$ is strongly
  connected and weight-balanced, we apply~\eqref{eq:LapBound} and the
  fact that $\ones_n^\top \xi_1 = 0$ to bound the above expression as
  \begin{multline*}
    - \frac{1}{2}\lambda_2(\Lap + \Lap^\top) \norm{\eta}^2 + \nu_1
    \eta^\top \xi_1 - \frac{1}{2} \beta \nu_1 \nu_2 \lambda_2(\Lap +
    \Lap^\top) \norm{\xi_1}^2 \notag
    \\
    - \alpha \norm{\xi_2}^2 - \nu_2 \xi_2^\top L \eta = \gamma^\top M
    \gamma ,
  \end{multline*}
  where $\eta = \zeta - \frac{1}{n} (\ones_n^\top \zeta)\ones_n$,
  $\gamma^\top = [\eta^\top, \xi_1^\top, \xi_2^\top]$, and
  \begin{align*}
    M =
    \begin{bmatrix}
      -\frac{1}{2} \lambda_2(\Lap + \Lap^\top) I & \frac{1}{2} \nu_1 I &
      - \frac{1}{2} \nu_2 \Lap^\top
      \\
      \frac{1}{2} \nu_1 I & -\frac{1}{2} \beta \nu_1 \nu_2 \lambda_2(\Lap
      + \Lap^\top) I & 0
      \\
      - \frac{1}{2} \nu_2 \Lap & 0 & -\alpha I
    \end{bmatrix}.
  \end{align*}
  Reasoning with the Schur complement~\citep{SB-LV:09}, $M \in
  \real^{3n \times 3n}$ is negative definite if
  \begin{align*}
    & -\frac{1}{2} \lambda_2(\Lap + \Lap^\top) I -
    \\
    & 
    \begin{bmatrix}
      \frac{1}{2} \nu_1 I
      &
      - \frac{1}{2} \nu_2 \Lap^\top
    \end{bmatrix}
    \begin{bmatrix}
      -\frac{1}{2} \beta \nu_1 \nu_2 \lambda_2(\Lap + \Lap^\top) I & 0
      \\
      0 & -\alpha
      I
    \end{bmatrix}^{-1}
    \begin{bmatrix}
      \frac{1}{2} \nu_1 I
      \\
      -\frac{1}{2} \nu_2 \Lap 
    \end{bmatrix}
    \\
    & = -\frac{1}{2} \lambda_2(\Lap + \Lap^\top) I + \frac{\nu_1}{2
       \beta \nu_2 \lambda_2(\Lap + \Lap^\top)} I +
       \frac{\nu_2^2}{4\alpha}
    \Lap^\top \Lap ,
  \end{align*}
  is negative definite. This latter fact is implied
  by~\eqref{eq:alpha-beta-cond-n}.  As a consequence, $
  \hat{\zeta}^\top g(P,\xi_1,\xi_2,\hat{\zeta}) \le 0$ and so,
  Lemma~\ref{le:continuity-set-valued-gen}\ref{as:continuity-1} holds.
  Moreover, $ \hat{\zeta}^\top g(P,\xi_1,\xi_2,\hat{\zeta}) = 0$ if
  and only if $\eta = \xi_1 = \xi_2 =0$, which means $\zeta \in
  \spn\{\ones_n\}$. Using this fact along with the definition of the
  set-valued Lie derivative and the characterization of
  optimizers~\eqref{eq:modified-edp-sol}, we deduce that $
  \hat{\zeta}^\top g(P,\xi_1,\xi_2,\hat{\zeta}) =0$ if and only if (a)
  $0 \in \SetLie_{\rLapgradnon} V_3(P,\xi_1,\xi_2)$ and (b) $P$ is a
  solution of the \ED problem. Fact (a) implies that
  Lemma~\ref{le:continuity-set-valued-gen}\ref{as:continuity-2} holds
  and hence,
  Proposition~\ref{pr:refined-lasalle-nonsmooth}\ref{as:refined-lasalle-2}
  holds too. Fact (b) implies that over the set $\HH_{P_l} \times
  \HH_0 \times \HH_0$, we have $0 \in \SetLie_{\rLapgradnon}
  V_3(P,\xi_1,\xi_2)$ if and only if $(P,\xi_1,\xi_2) \in
  \bFFa$. Since, $\bFFa$ belongs to a level set of~$V_3$, we conclude
  that
  Proposition~\ref{pr:refined-lasalle-nonsmooth}\ref{as:refined-lasalle-1}
  holds too.
 
  To be able to apply Proposition~\ref{pr:refined-lasalle-nonsmooth}
  and conclude the proof, it remains to show that the trajectories
  starting from $D(\real^n \times \real^n \times \HH_0)$ are bounded.
  We reason by contradiction, i.e., assume there exists a trajectory
  $t \mapsto (P(t),\xi_1(t),\xi_2(t))$, with initial condition
  $(P(0),\xi_1(0),\xi_2(0)) \in D(\real^n \times \real^n \times
  \HH_0)$
  of $\rLapgradnon$ such that $\norm{(P(t),\xi_1(t),\xi_2(t)} \to
  \infty$.  Since $V_3$ is radially unbounded, this implies
  $V_3(P(t),\xi_1(t),\xi_2(t)) \to \infty$.  Additionally, from
  Lemma~\ref{le:omega-limit-set-edp}, we know that $\ones_n^\top P(t)
  \to P_l$ and $\ones_n^\top \xi_1(t) \to 0$.  Thus, there exists a
  sequence of times $\{t_k\}_{k=1}^{\infty}$ with $t_k \to \infty$
  such that for all $k \in \integerspositive$,
  \begin{subequations}\label{eq:xi-bound}
    \begin{align}
      \abs{\ones_n^\top \xi_1(t_k)} & < {1}/{k}
      , \label{eq:xi-bound-1}
      \\
      \max \SetLie_{\rLapgradnon} V_3(P(t_k),\xi_1(t_k),\xi_2(t_k)) & >
      0. \label{eq:xi-bound-2}
    \end{align}
  \end{subequations}
  Note that~\eqref{eq:xi-bound-2} implies that there exists a sequence
  $\{\zeta_k\}_{k = 1}^{\infty}$ with $\zeta_k \in \partial
  f^{\eps}(P(t_k))$ such that
  \begin{multline*}
    -\zeta_k^\top \Lap \zeta_k + \nu_1 \zeta_k^\top \xi_1(t_k) - \beta
    \nu_1 \nu_2 \xi_1(t_k)^\top \Lap \xi_1(t_k)
    \\
    - \alpha \norm{\xi_1(t_k)}^2 - \nu_2 \xi_2(t_k)^\top \Lap \zeta_k
    > 0 ,
  \end{multline*}
  for all $k \in \integerspositive$, where we have used the fact that
  an element of $\SetLie_{\rLapgradnon} V_3 (P,\xi_1,\xi_2)$ has the
  form given in~\eqref{eq:liederivative-element}.  Letting $\eta_k =
  \zeta_k - \frac{1}{n}(\ones_n^\top \zeta_k) \ones_n$, we
  use~\eqref{eq:LapBound} to deduce from the above inequality that
  \begin{align*}
    & -\frac{1}{2} \lambda_2(\Lap + \Lap^\top) \norm{\eta_k}^2 + \nu_1
    \eta_k^\top \xi_1(t_k) + \frac{1}{n} \nu_1 (\ones_n^\top \zeta_k)
    (\ones_n^\top \xi_1(t_k))
    \\
    & - \frac{1}{2} \beta \nu_1 \nu_2 \lambda_2(\Lap + \Lap^\top)
    \norm{\xi_1(t_k) - \frac{1}{n} (\ones_n^\top \xi_1(t_k))
      \ones_n}^2
    \\
    & - \alpha \norm{\xi_1(t_k)}^2 - \nu_2 \xi_2(t_k)^\top \Lap
    \eta_k > 0.
  \end{align*}
  Further, using the expression
  \begin{align*}
    \norm{\xi_1(t_k) - \frac{1}{n} (\ones_n^\top \xi_1(t_k))
      \ones_n}^2 = \norm{\xi_1(t_k)}^2 - \frac{1}{n} (\ones_n^\top
    \xi_1(t_k))^2,
  \end{align*}
  the inequality can be rewritten as
  \begin{multline*}
    \gamma_k^\top M \gamma_k + \frac{1}{n} \nu_1 (\ones_n^\top \zeta_k)
    (\ones_n^\top \xi_1(t_k))
    \\
    + \frac{\beta \nu_1 \nu_2}{2 n } \lambda_2(\Lap + \Lap^\top)
    (\ones_n^\top \xi_1(t_k))^2 > 0 ,
  \end{multline*}
  where $\gamma_k^\top = [ \eta_k^\top, \, \xi_1(t_k)^\top, \,
  \xi_2(t_k)^\top ]$.  Using now the bound~\eqref{eq:xi-bound-1}, we
  arrive at the inequality,
  \begin{equation}\label{eq:lie-bound}
    \gamma_k^\top M \gamma_k + \frac{\nu_1}{nk} \abs{\ones_n^\top \zeta_k} +
    \frac{\beta \nu_1 \nu_2}{2 n k^2} \lambda_2(\Lap + \Lap^\top) > 0.
  \end{equation}
  Next, we consider two cases, depending on whether the sequence
  $\{P(t_k)\}$ is (a) bounded or (b) unbounded. In case (a), the
  sequence $\{(\xi_1(t_k),\xi_2(t_k))\}$ must be unbounded. Since $M$
  is negative definite, we have $\gamma_k^\top M \gamma_k \le
  \lambda_{\max}(M) \norm{(\xi_1(t_k), \xi_2(t_k))}^2$.
  Thus,~\eqref{eq:lie-bound} implies that
  \begin{align*}
    \lambda_{\max} (M) \norm{(\xi_1(t_k), \xi_2(t_k))}^2 & +
    \frac{\nu_1}{nk} \abs{\ones_n^\top \zeta_k}
    \\
    & + \frac{\beta \nu_1 \nu_2}{2 n k^2} \lambda_2(\Lap + \Lap^\top) >
    0.
  \end{align*} 
  Now, from the expression of $\partial f^\eps$, since $\{P(t_k)\}$ is
  bounded, the sequence $\{\zeta_k\}$ must be bounded. Combining these
  facts with $\lambda_{\max} (M) <0$, one can find $\bar{k} \in
  \integerspositive$ such that the above inequality is violated for
  all $k \ge \bar{k}$, which is a contradiction.  For case (b), we use
  the bound $\gamma_k^\top M \gamma_k \le
  \lambda_{\max}(M)\norm{\eta_k}^2$ to deduce
  from~\eqref{eq:lie-bound} that
  \begin{align*}
    \lambda_{\max} (M) \norm{\eta_k}^2 + \frac{\nu_1}{nk}
    \abs{\ones_n^\top \zeta_k} + \frac{\beta \nu_1 \nu_2}{2 n k^2}
    \lambda_2(\Lap + \Lap^\top) > 0.
  \end{align*}
  One can then use a similar argument as laid out in the
  proof of Theorem~\ref{th:central-convergence}, considering the two
  cases of $\abs{\ones_n^\top \zeta_k}$ being bounded or unbounded,
  arriving in both cases at similar contradictions. This concludes the
  proof. \qed
\end{pf}

Note that as a consequence of the above result, the \dacplap dynamics
does not require any specific pre-processing for the initialization of
the power allocations. Each generator can select any generation level,
independent of the other units, and the algorithm guarantees
convergence to the solutions of the \ED problem.

\begin{remark}\longthmtitle{Distributed selection of algorithm design
    parameters}\label{re:dist-select}
  {\rm The convergence of the \dacplap dynamics relies on a selection
    of the parameters $\alpha$, $\beta$, $\nu_1$ and $\nu_2 \in
    \realpositive$ that satisfy~\eqref{eq:alpha-beta-cond-n}. Checking
    this inequality requires knowledge of the spectrum of matrices
    related to the Laplacian matrix, and hence the entire network
    structure.  Here, we provide an alternative condition that
    implies~\eqref{eq:alpha-beta-cond-n} and can be checked by the
    units in a distributed way.  Let $\nm$ be an upper bound on the
    number of units, $ \wdoutmax$ be an upper bound on the out-degree
    of all units, and $\amin$ be a lower bound on the edge weights,
    \begin{align}\label{eq:req-bounds}
      n \le \nm, \; \max_{i \in \vertices} \wdout(i) \le \wdoutmax, \;
      \min_{(i,j) \in \edges} a_{ij} \ge \amin.
    \end{align}
    A straightforward generalization of~\citep[Theorem 4.2]{BM:91-gc}
    for weighted graphs gives rise to the following lower bound on
    $\lambda_2(\Lap+\Lap^\top)$,
    \begin{equation}\label{eq:graph-bounds1}
      \frac{4 \amin}{\nm^2}  \le \lambda_2(\Lap +
      \Lap^\top).
    \end{equation}
    On the other hand, using properties of matrix 
    norms~\citep[Chapter 9]{DSB:05}, one can deduce
    \begin{align}
      \lambda_{\max}(\Lap^\top \Lap) & = \norm{\Lap}^2  \le (\sqrt{n}
      \norm{\Lap}_{\infty})^2 \notag
      \\
      & \le (2 \sqrt{n} \wdoutmax)^2 \le 4 \nm (\wdoutmax)^2.
      \label{eq:graph-bounds2}
    \end{align}
    Using~\eqref{eq:graph-bounds1}-\eqref{eq:graph-bounds2}, the
    left-hand side of~\eqref{eq:alpha-beta-cond-n} can be upper bounded
    by
    \begin{align*}
      \frac{\nu_1}{\beta \nu_2 \lambda_2(\Lap + \Lap^\top)} & +
      \frac{\nu_2^2 \lambda_{\max}(L^\top L)}{2 \alpha}
      \\
      & \qquad \le \frac{\nu_1 \nm^2}{4 \amin \beta \nu_2} + \frac{2
        \nu_2^2 \nm (\wdoutmax)^2}{\alpha}.
    \end{align*}
    Further, the right-hand side of~\eqref{eq:alpha-beta-cond-n} can be
    lower bounded using~\eqref{eq:graph-bounds1}. Putting the two
    together, we obtain the new condition
    \begin{equation}\label{eq:alpha-beta-cond-d}
      \frac{\nu_1 \nm^2}{4 \amin \beta \nu_2} + \frac{2 \nu_2^2 \nm
        (\wdoutmax)^2}{\alpha} < 
      \frac{4 \amin}{\nm^2} ,
    \end{equation}
    which implies~\eqref{eq:alpha-beta-cond-n}.  The network can
    ensure that this condition is met in various ways.  For instance,
    if the bounds $\nm$, $ \wdoutmax$, and $\amin$ are not available,
    the network can implement distributed algorithms for max- and
    min-consensus~\citep{WR-RWB:08} to compute them in finite time.
    Once known, any generator can select $\alpha$, $\beta$, $\nu_1$
    and $\nu_2$ satisfying~\eqref{eq:alpha-beta-cond-d} and broadcast
    its choice.  Alternatively, the computation of the design
    parameters can be implemented concurrently with the determination
    of the bounds via consensus by specifying a specific formula to
    select them that is guaranteed to
    satisfy~\eqref{eq:alpha-beta-cond-d}.  Note that the units
    necessarily need to agree on the parameters, otherwise if each
    unit selects a different set of parameters, the dynamic average
    consensus would not track the average input signal.  } \oprocend
\end{remark}

\begin{remark}\longthmtitle{Distributed loads and
    transmission losses}\label{re:load-agg-2}
  {\rm Here we expand on our observations in
    Remark~\ref{re:load-agg-1} regarding the inclusion of additional
    constraints on the \ED problem.  Our algorithmic solution can be
    easily modified to deal with the alternative scenarios studied
    in~\citep{ZZ-XY-MC:11,SK-GH:12,GB-AD-FLL-DN-BT:14,VL-AV:13}, where
    each generator has the knowledge of the load at the corresponding
    bus that it is connected to and the total load is the aggregate of
    these individual loads.  Mathematically, denoting the load
    demanded at generator bus $i$ by $P^L_i \in \real$, the total load
    is given by $P_l = \sum_{i=1}^n P^L_i$. For this case, replacing
    the vector $P_l e_r$ by $P^L$ in the \dacplap
    dynamics~\eqref{eq:red-two} gives an algorithm that solves the \ED
    problem for the load~$P_l$.  Our solution strategy can also handle
    transmission losses as modeled in~\citep{GB-AD-FLL-DN-BT:14},
    where it is assumed that each generator~$i$ can estimate the power
    loss in the transmission lines adjacent to it. With those values
    available, the generator could add them to the quantity $P^L_i$,
    which would make the network find a power allocation that takes
    care of the transmission losses.  } \oprocend
\end{remark}

\subsection{Robustness analysis}\label{subsec:robust-prop}

In this section, we study the robustness properties of the \dacplap
dynamics in the presence of time-varying load signals and intermittent
power unit generation.  Our analysis relies on the exponential
stability of the mismatch dynamics between total generation and load
established in Lemma~\ref{le:omega-limit-set-edp}, which implies
that~\eqref{eq:mismatch-dyn} is input-to-state stable
(ISS)~\citep[Lemma 4.6]{HKK:02}, and consequently robust against
arbitrary bounded perturbations.  The following result provides an
explicit, exponentially decaying, bound for the evolution of any
trajectory of~\eqref{eq:mismatch-dyn}.
While the rate of decay can also be determined by computing the
eigenvalues of matrix defining the dynamics, here we employ a Lyapunov
argument to obtain also the value of the gain associated to the rate.

\begin{lemma}\longthmtitle{Convergence rate of the mismatch
    dynamics~\eqref{eq:mismatch-dyn}}\label{le:mismatch-exp-stable1}
  Let $R \in \real^{2 \times 2}$ be defined by
  \begin{align*}
    R = \frac{1}{2\alpha \nu_1 \nu_2} 
    \begin{bmatrix}
      \alpha^2+\nu_1 \nu_2 + (\nu_1 \nu_2)^2 & \alpha
      \\
      \alpha & 1+\nu_1 \nu_2
    \end{bmatrix} .
  \end{align*}
  Then $R \succ 0$ and any trajectory $t \mapsto x(t)$ of the
  dynamics~\eqref{eq:mismatch-dyn} satisfies $ \norm{x(t)} \le c_1
  e^{-c_2t} \norm{x(0)}$, where $c_1 =
  \sqrt{{\lambda_{\max}(R)}/{\lambda_{\min}(R)}}$ and $c_2 =
  {1}/{2\lambda_{\max}(R)}$.
\end{lemma}
\begin{pf}
  Let $A \in \real^{2\times 2}$ be the system matrix
  of~\eqref{eq:mismatch-dyn}. Then, one can see that $A^\top R + R A =
  -I$, i.e., $V_4(x) = x^\top R x$ is a Lyapunov function
  for~\eqref{eq:mismatch-dyn}.  Note that
  \begin{equation}\label{eq:lyap-norm-bounds}
    \lambda_{\min}(R) \norm{x}^2 \le V_4(x) \le \lambda_{\max}(R)
    \norm{x}^2.
  \end{equation} 
  From the Lyapunov equation, we have $\Lie_{Ax} V_4(x) = - \norm{x}^2
  \le - \frac{1}{\lambda_{\max}(R)} V_4(x)$, which implies $V_4(x(t)) \le
  e^{-1/\lambda_{\max}(R)} V_4(x(0))$ along any trajectory $t \mapsto
  x(t)$ of~\eqref{eq:mismatch-dyn}. Again
  using~\eqref{eq:lyap-norm-bounds}, we get
  \begin{align*}
    \norm{x(t)}^2 \le \frac{\lambda_{\max}(R)}{\lambda_{\min}(R)}
    e^{-1/\lambda_{\max}(R)} \norm{x(0)}^2,
  \end{align*}
  which concludes the claim. \qed
\end{pf}


In the above result, it is interesting to note that the convergence
rate is independent of the specific communication digraph (as long as
it is weight-balanced). We use next the exponentially decaying bound
obtained above to illustrate the extent to which the network can
collectively track a dynamic load (which corresponds to a time-varying
perturbation in the mismatch dynamics) and is robust to 
intermittent power generation (which corresponds to
perturbations in the state of the mismatch dynamics).

\subsubsection{Tracking dynamic loads}\label{subssec:tracking}

Here we consider a time-varying total load given by a twice
continuously differentiable trajectory $\realnonnegative \ni t \mapsto
P_l(t)$ and show how the total generation of the network under the
\dacplap dynamics tracks~it. We assume the signal is known to an
arbitrary unit $r \in \until{n}$. In this case, the
dynamics~\eqref{eq:mismatch-dyn} takes the following form
\begin{align*}
  \begin{bmatrix}
    \dot x_1
    \\
    \dot x_2
  \end{bmatrix}
  = 
  \begin{bmatrix}
    0 & 1 
    \\ 
    - \nu_1 \nu_2 & - \alpha
  \end{bmatrix}
  \begin{bmatrix}
    x_1
    \\
    x_2
  \end{bmatrix}
  + 
  \begin{bmatrix}
    0
    \\
    -\alpha \dot P_l - \ddot P_l
  \end{bmatrix}.
\end{align*}
Using Lemma~\ref{le:mismatch-exp-stable1}, one can compute the
following bound on any trajectory of the above system 
\begin{align*}
  \norm{x(t)} \le c_1 e^{-c_2 t} \norm{x(0)} + \frac{c_1}{c_2} \sup_{s
    \in [0,t]} \abs{\alpha \dot P_l(s) + \ddot P_l(s)}.
\end{align*}
In particular, for a signal with bounded $\dot P_l$ and $\ddot P_l$,
the mismatch between generation and load, i.e., $x_1(t)$ is
bounded. Also, the mismatch has an ultimate bound as $t \to
\infty$. The following result summarizes this notion formally. The
proof is straightforward application of
Lemma~\ref{le:mismatch-exp-stable1} following the exposition of
input-to-state stability in~\citep{HKK:02}.

\begin{proposition}\longthmtitle{Power mismatch is ultimately bounded
    for dynamic load under \dacplap
    dynamics}\label{pr:mismatch-ult-bound}
  Let $\realnonnegative \ni t \mapsto P_l(t)$ be twice continuously
  differentiable and such that
  \begin{align*}
    \sup_{t \ge 0} \abs{\dot P_l(t)} \le d_1, \quad \sup_{t \ge 0}
    \abs{\ddot P_l(t)} \le d_2 ,
  \end{align*}
  for some $d_1 , d_2 > 0$. Then, the mismatch $\ones_n^\top P(t) -
  P_l(t)$ between load and generation is bounded along the
  trajectories of~\eqref{eq:robust-edp-dist} and has ultimate bound
  $\frac{c_1}{c_2} (\alpha d_1 + d_2)$, with $c_1$, $c_2$ given in
  Lemma~\ref{le:mismatch-exp-stable1}. Moreover, if $\dot P_l(t) \to
  0$ and $\ddot P_l(t) \to 0$ as $t \to \infty$, then $\ones_n^\top
  P(t) \to P_l(t)$ as $t \to \infty$.
\end{proposition}

\subsubsection{Robustness to intermittent power
  generation}\label{subssec:robust-add-delete}

Here, we characterize the algorithm robustness against unit addition
and deletion to capture scenarios with intermittent power generation.
Addition and deletion events are modeled via a time-varying
communication digraph, which we assume remains strongly connected and
weight-balanced at all times. When a unit stops generating power
(deletion event), the corresponding vertex and its adjacent edges are
removed.  When a unit starts providing power (addition event), the
corresponding node is added to the digraph along with a set of edges.
Given the intricacies of the convergence analysis for the \dacplap
dynamics, cf. Theorem~\ref{th:convergence-rlngd}, it is important to
make sure that the state $v$ remains in the set $\HH_0$,
irrespectively of the discontinuities caused by the events.  The
following routine makes sure that this is the case.
\begin{quote}
  \tinv: When a unit $i$ joins the network at time $t$, it starts with
  $v_i(t) = 0$.
  When a unit $i$ leaves the network at time $t$, it passes a token
  with value $v_i(t)$ to one of its in-neighbors $j \in \Nin(i)$, who
  resets its value to $v_j(t) + v_i(t)$.
\end{quote}
The \tinv routine ensures that the dynamics~\eqref{eq:mismatch-dyn} is
the appropriate description for the evolution of the load satisfaction
mismatch. This, together with the ISS property established in
Lemma~\ref{le:mismatch-exp-stable1}, implies that the mismatch effect
in power generation caused by addition/deletion events vanishes
exponentially fast. In particular, if the number of addition/deletion
events is finite, then the set of generators converge to the solution
of the \ED problem.  We formalize this next.

\begin{proposition}\longthmtitle{Convergence of
    \dacplap dynamics under intermittent
    power generation}\label{pr:intermittent-gen}
  Let $\nm$ be the maximum number of generators that can contribute to
  the power generation at any time. Let $\Sigma_{\nm}$ be the set of
  digraphs that are strongly connected and weight-balanced and whose
  vertex set is included in $\until{\nm}$. Let $\sigma:[0,\infty) \to
  \Sigma_{\nm}$ be a piecewise constant, right-continuous switching
  signal described by the set of switching times $ \{t_1, t_2, \dots
  \} \subset \realnonnegative$, with $t_k \le t_{k+1}$, each
  corresponding to either an addition or a deletion event.  Denote by
  $\rLapgradnon^{\sigma}$ the switching \dacplap dynamics
  corresponding to $\sigma$, defined by~\eqref{eq:robust-edp-dist}
  with $\Lap$ replaced by $\Lap(\sigma(t))$ for all $t \ge 0$, and
  assume agents execute the \tinv routine when they leave or join the
  network. Then,
  \begin{enumerate}
  \item at any time $ t \in \{0\} \cup \{t_1, t_2, \dots \}$, if the
    variables $(P(t),z(t))$ for the generators in $\sigma(t)$ satisfy
    $\abs{\ones_n^\top P(t) - P_l} \le M_1$ and $\abs{\ones_n^\top
      z(t)} \le M_2$ for some $M_1, M_2 > 0$, then the magnitude of
    the mismatch between generation and load becomes less than or
    equal to $\rho>0$ in time
    \begin{align*}
      t_{\rho} = \frac{1}{c_2}\ln \Bigl(\frac{c_1(M_1+\nu_1 M_2)}
      {\rho}\Bigr) ,
    \end{align*}
    provided no event occurs in the interval $(t,t+t_{\rho})$;
  \item if the number of events is finite, say $N$, then the
    trajectories of $\rLapgradnon^{\sigma}$ converge to the set of
    solutions of the \ED problem for the group of generators
    in~$\sigma(t_N)$ provided~\eqref{eq:alpha-beta-cond-n} is met for
    $\sigma(t_N)$.
  \end{enumerate}
\end{proposition}

Note that the generators can ensure that the
condition~\eqref{eq:alpha-beta-cond-n}, required for the convergence
of the \dacplap dynamics, holds at all times even under addition and
deletion events, if they rely on verifying
that~\eqref{eq:alpha-beta-cond-d} holds and the
bounds~\eqref{eq:req-bounds} are valid for all the topologies in
$\Sigma_{\nm}$.

\section{Simulations in a IEEE 118 bus system}\label{sec:sims}

\begin{figure*}[htb]
  \centering \subfloat[Power allocation]{\includegraphics[width = 0.3
    \linewidth]{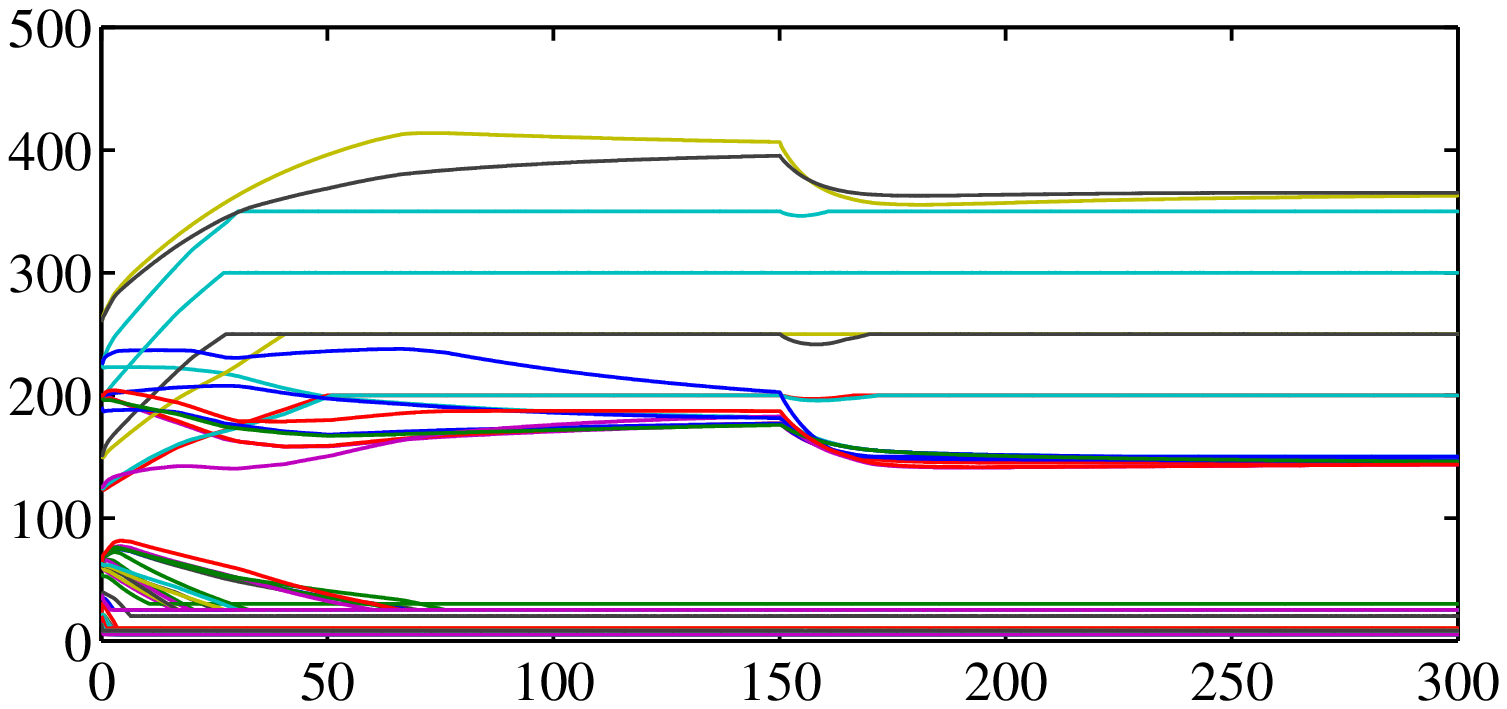}} \quad \subfloat[Total
  cost]{\includegraphics[width = 0.3
    \linewidth]{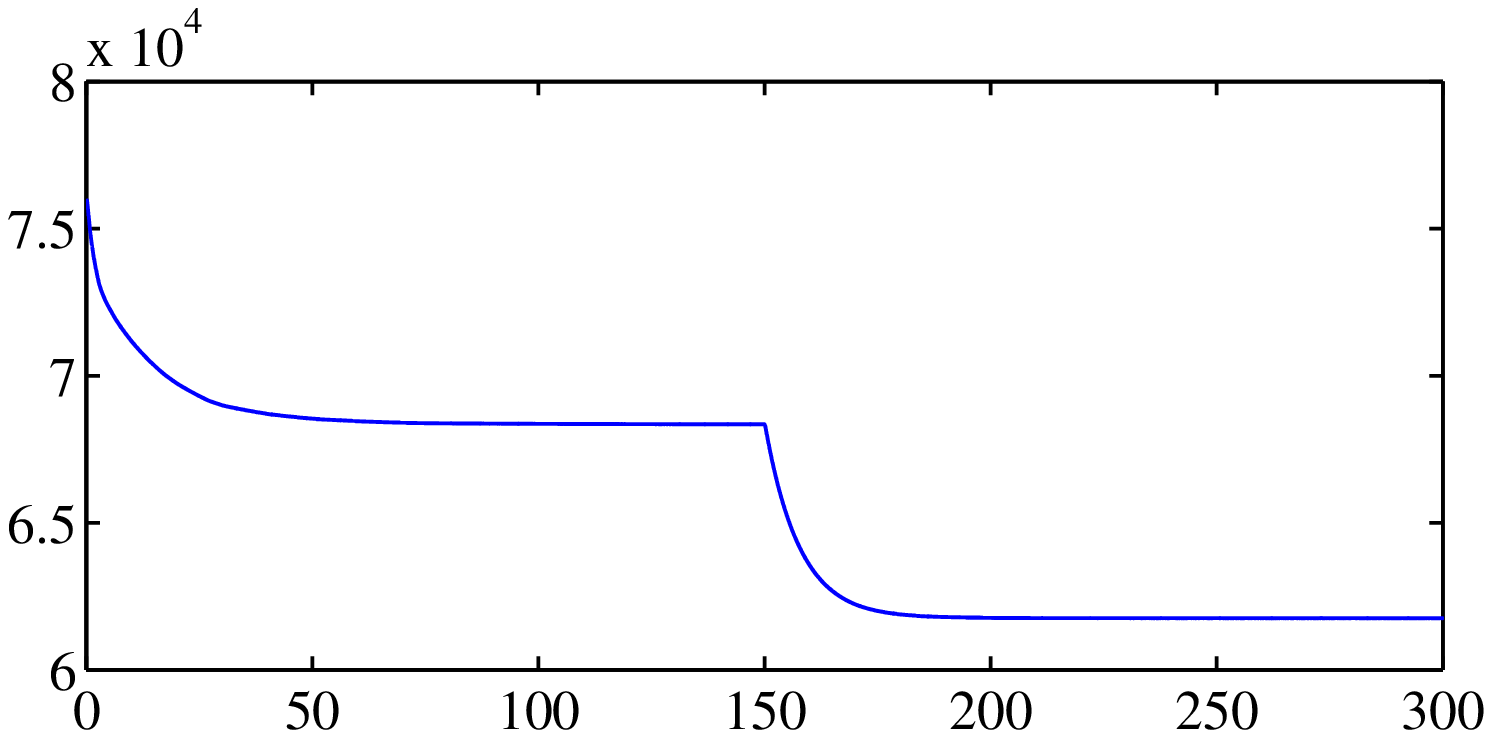}} \quad \subfloat[Total
  mismatch]{\includegraphics[width = 0.3
    \linewidth]{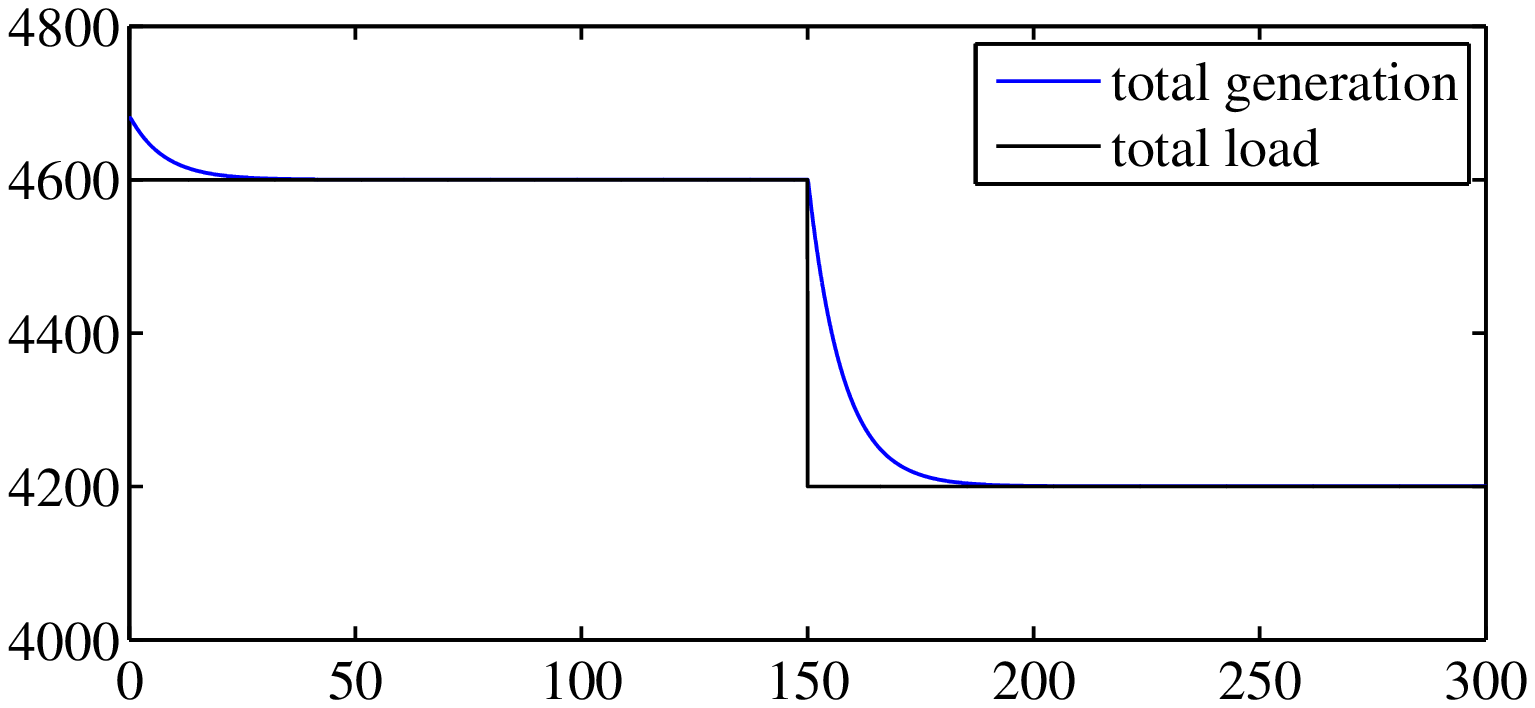}}
  \\
  \subfloat[Power allocation]{\includegraphics[width = 0.3
    \linewidth]{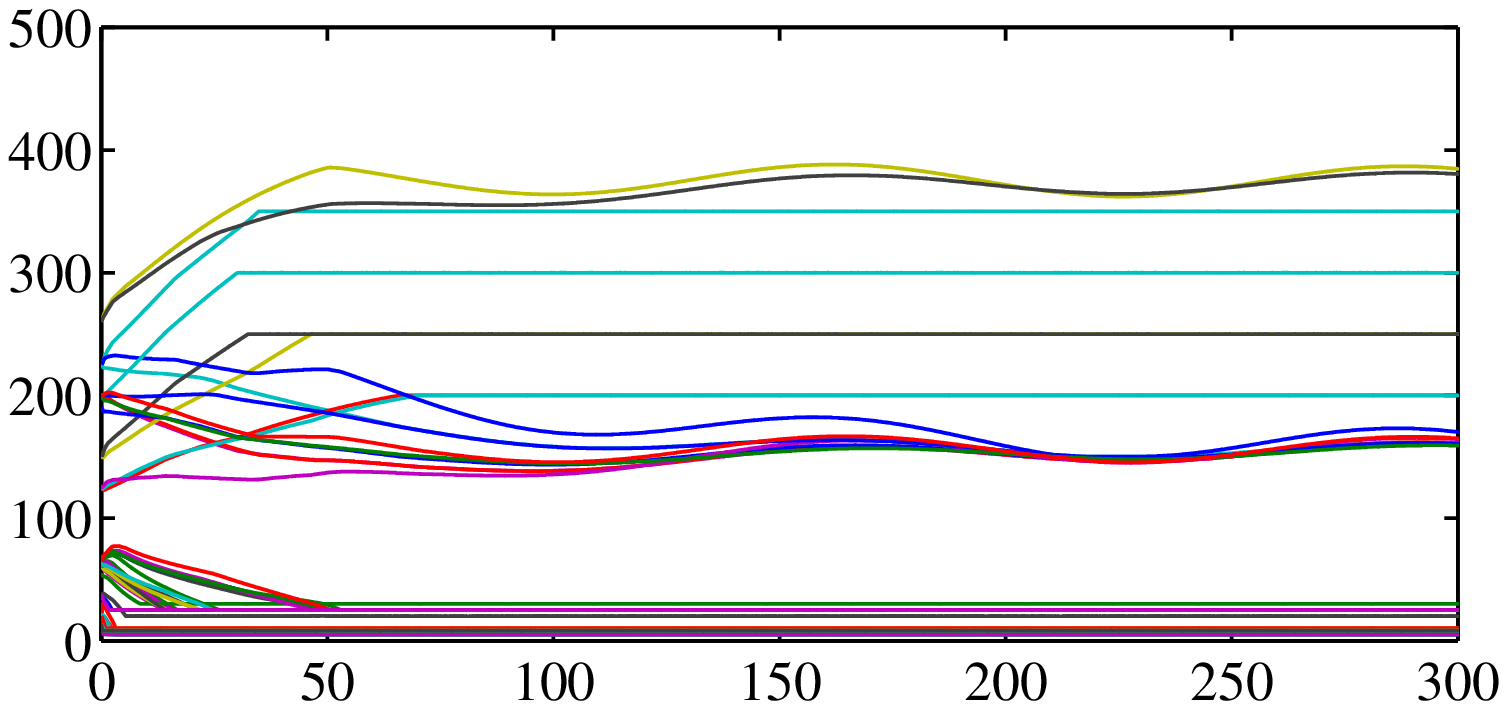}} \quad \subfloat[Total
  cost]{\includegraphics[width = 0.3
    \linewidth]{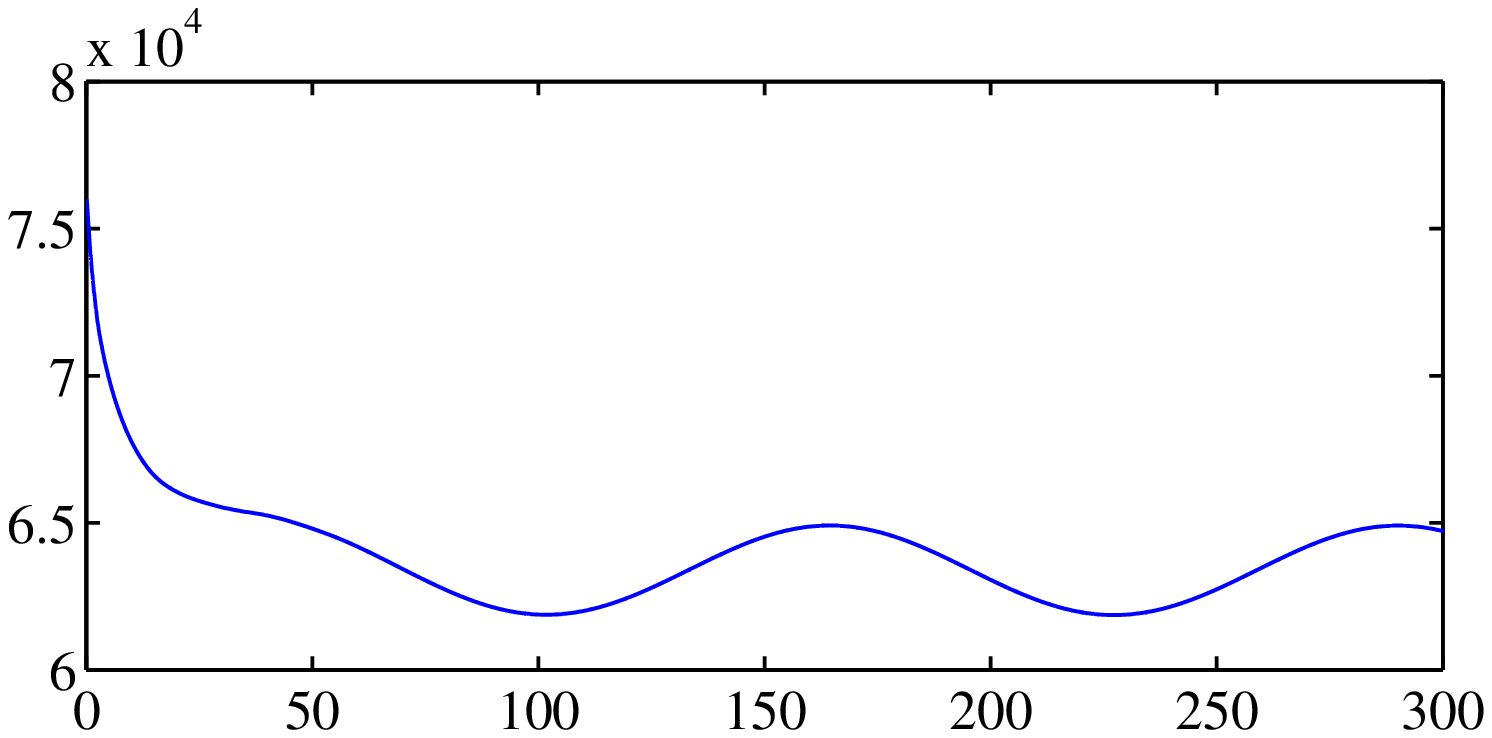}} \quad \subfloat[Total
  mismatch]{\includegraphics[width = 0.3
    \linewidth]{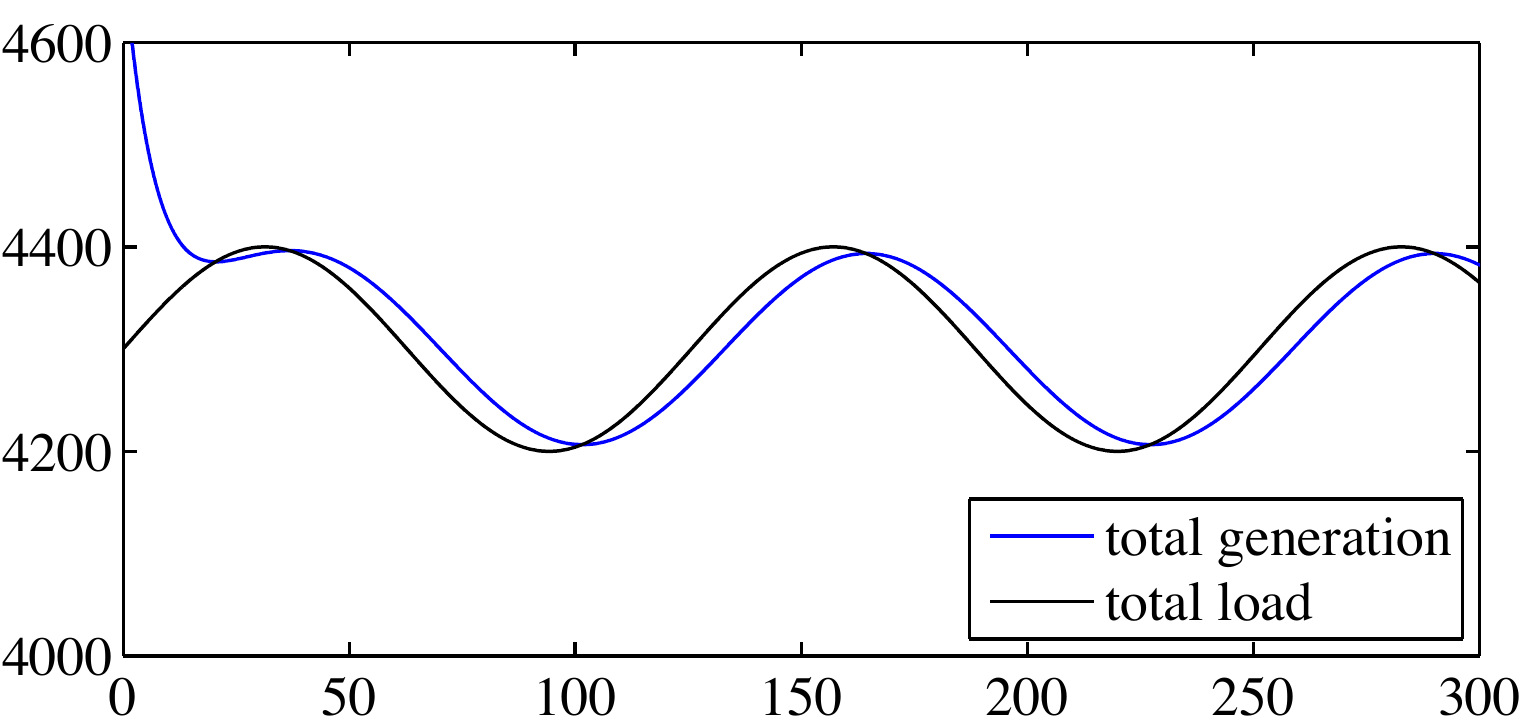}}
  \\
  \subfloat[Power allocation]{\includegraphics[width = 0.3
    \linewidth]{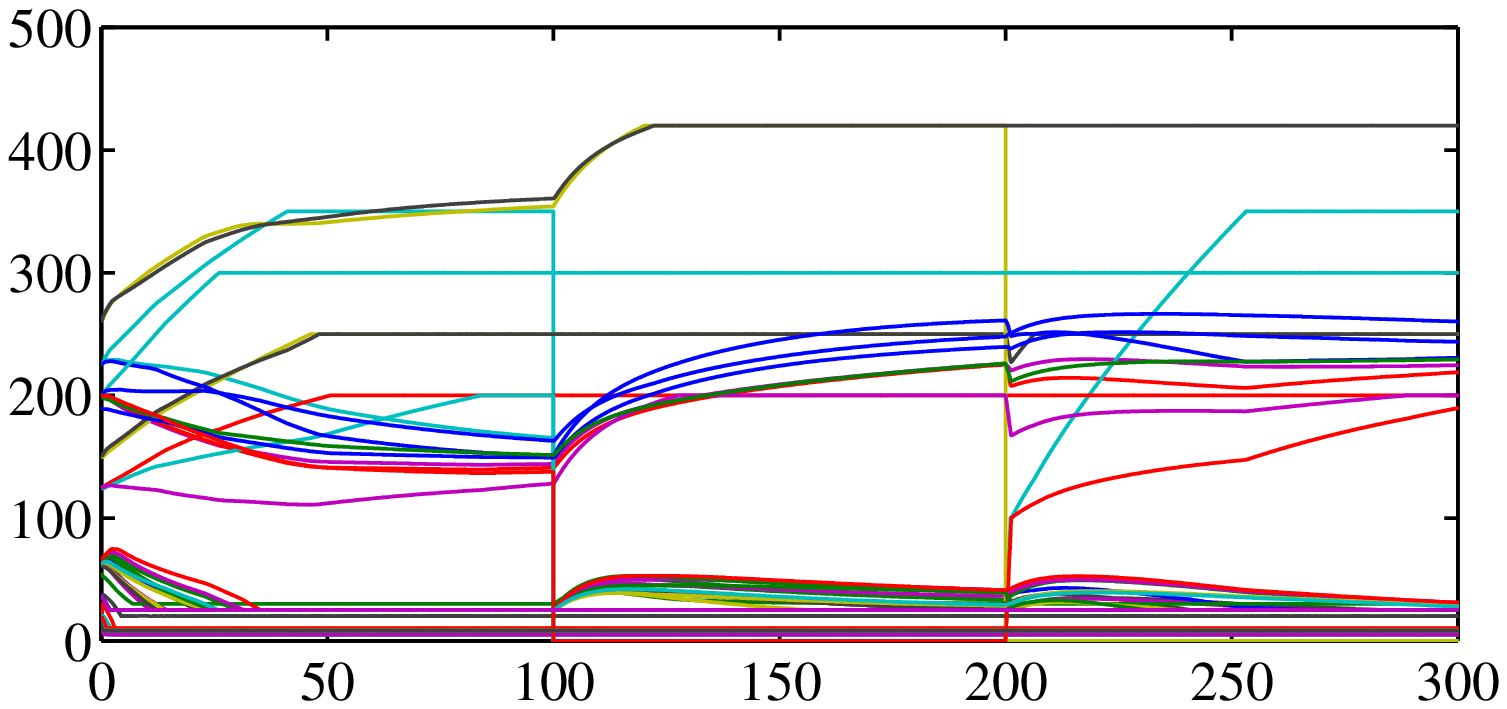}} \quad \subfloat[Total
  cost]{\includegraphics[width = 0.3
    \linewidth]{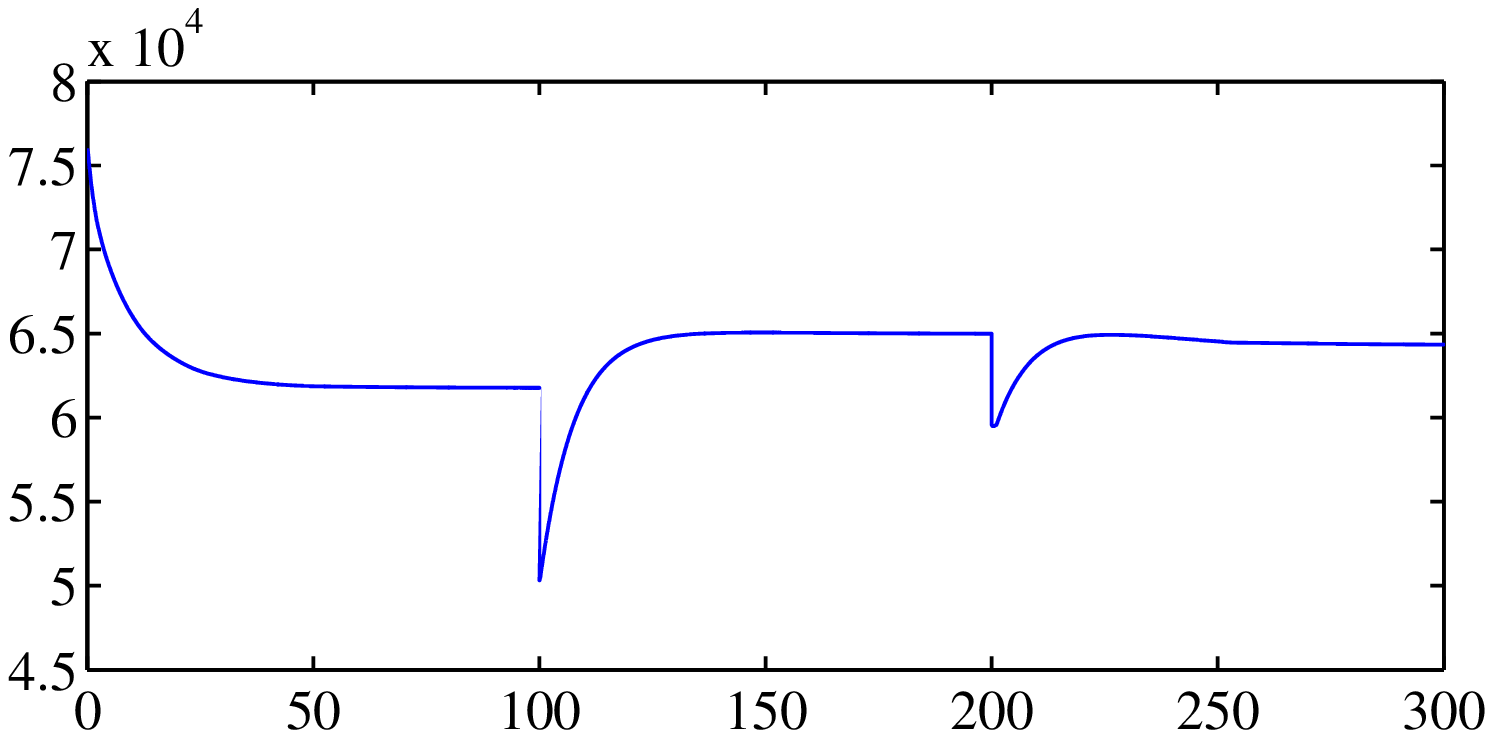}} \quad \subfloat[Total
  mismatch]{\includegraphics[width = 0.3
    \linewidth]{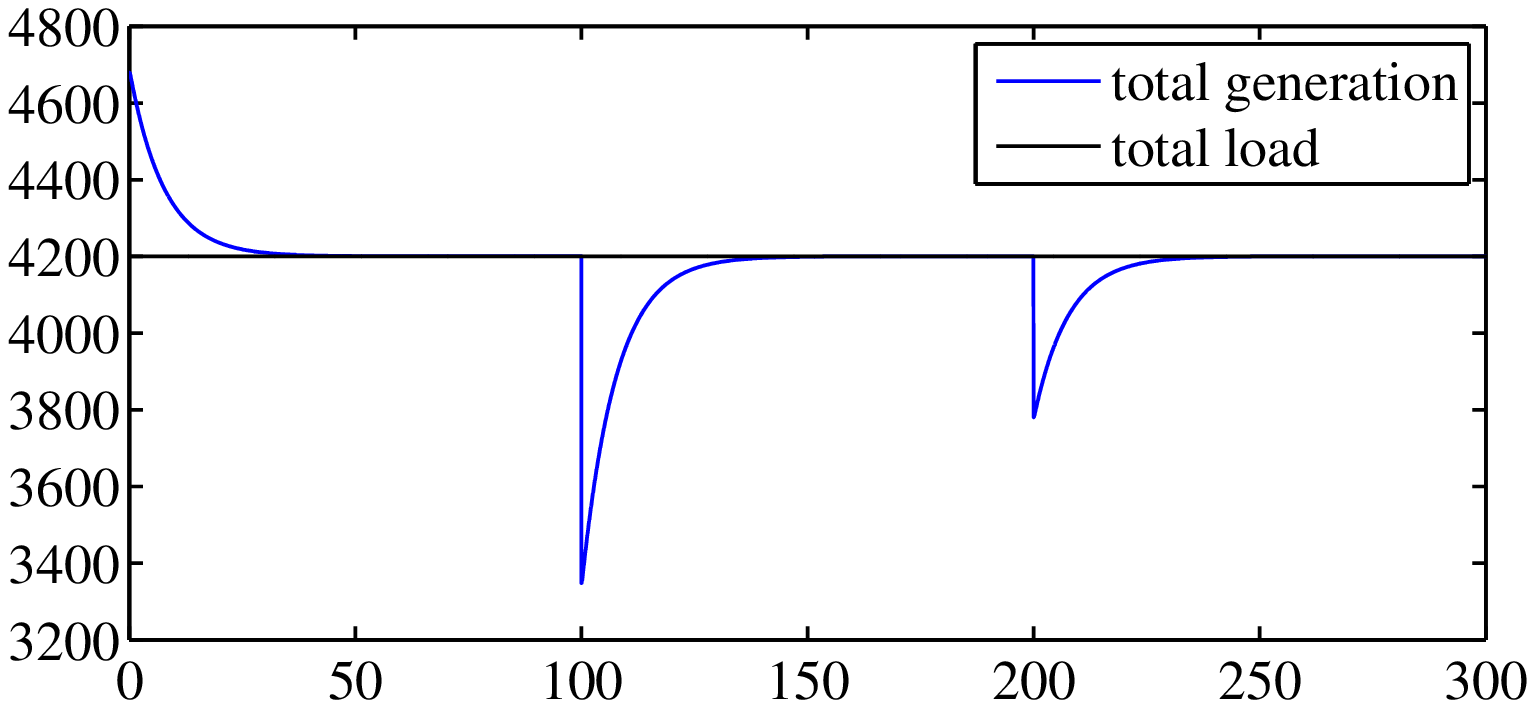}}
  \caption{Evolution of the power allocation, the total cost, and the
    total mismatch between generation and load under the \dacplap
    dynamics for the IEEE 118 bus example in different scenarios. In
    the first case (a)-(c), the communication topology is $\GGa$, the
    load is initially $4600$ and later $4200$, and parameters are
    $\nu_1 = 1$, $\nu_2 = 1.3$, $\alpha = 10$, $\beta = 40$, and $\eps
    = 0.0086$.  In the second scenario (d)-(f), the digraph and the
    parameters remain the same but the load is time-varying given by
    $P_l(t) = 4300 + 100 \sin(0.05t)$.  In the last case (g)-(i), the
    parameters remain the same, the communication graph is initially
    the graph $\GGau$. At $t = 100s$, units $\{4,11,25,45\}$ leave the
    network, resulting in the communication topology $\GGi$, and the
    remaining agents run the \tinv routine.  Later, at $t = 200s$,
    units $\{11,45\}$ join the network while unit $27$ leaves it,
    resulting in the communication topology $\GGf$. After implementing
    the \tinv routine, the \dacplap dynamics  eventually converges to
    an optimizer of the \ED problem for the network $\GGf$.
  }
  \label{fig:one}
\end{figure*}

This section illustrates the convergence of the \dacplap dynamics to
the solutions of the \ED problem~\eqref{eq:EDP} starting from any
initial power allocation and its robustness properties.  We consider
the IEEE 118 bus system~\citep{IEEE-118}, that consists of 54
generators. The cost function of each generator $i$ is quadratic,
$f_i(P_i) = a_i + b_i P_i + c_i P_i^2$, with coefficients belonging to
the ranges $a_i \in [6.78, 74.33]$, $b_i \in [8.3391,37.6968]$, and
$c_i \in [0.0024, 0.0697]$. The communication topology is the digraph
$\GGa$ described in Table~\ref{tb:graphs}.
\begin{table*}
  \centering
  {\small
    \begin{tabular}{|l|l|}
      \hline
      $\GGa$ & digraph over $54$ vertices consisting of a directed cycle
      through vertices $1,\dots, 54$ and bi-directional edges 
      \\[-1ex]
      & $\{(i,\idf(i+5)), (i,\idf(i+10)), (i,\idf(i+15)),
      (i,\idf(i+20))\}$ for each 
      $i \in \until{54}$, where
      \\[-1ex]
      &   $\idf(x) = x$ if $x \in \until{54}$ and
      $ x-54$ otherwise.   All edge weights  are 0.1.   
      \\
      \hline
      $\GGau$ & obtained from $\GGa$ by
      replacing the directed cycle with  an undirected one keeping the 
      edge weights same \\
      \hline
      $\GGi$ & obtained from $\GGau$ by removing
      the  vertices $\{4,11,25,45\}$ and the edges adjacent to them 
      \\
      \hline
      $\GGf$ & obtained from $\GGau$ by 
      removing the vertices $\{4,25,27\}$ and the edges adjacent to them
      \\
      \hline
    \end{tabular}
  }
  \caption{Definition of the digraphs $\GGa$, $\GGau$, $\GGi$, and
    $\GGf$.}\label{tb:graphs}
\end{table*}
We choose the design parameters as $\nu_1 = 1, \nu_2 = 1.3, \alpha =
10, \beta = 40, \eps = 0.0086$, which satisfy the
conditions~\eqref{eq:EpsBound} and~\eqref{eq:alpha-beta-cond-n}
for~$\GGa$. The total load is $4600$ for the first $150$ seconds and
$4200$ for the next $150$ seconds, and is known to unit~$3$.
Figure~\ref{fig:one}(a)-(c) depicts the evolution of the power
allocation, total cost, and the mismatch between the total generation
and load under the \dacplap dynamics starting at the initial condition
$(P(0),z(0),v(0)) = (0.5*(P^m+P^M),0,0)$.  Note that the generators
initially converge to a power allocation that meets the load $4600$
and minimizes the total cost of generation. Later, with the decrease
in desired load to $4200$, the network decreases the total generation
while minimizing the total cost.

\begin{figure}
  \centering
  \includegraphics[width = 0.7
  \linewidth]{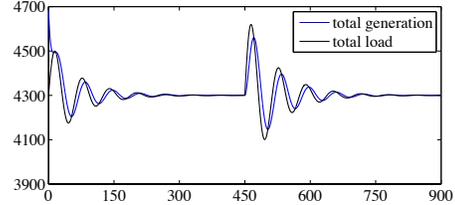}
  \caption{Evolution of the total power generation for the IEEE 118
    bus example under the \dacplap dynamics for the communication
    digraph $\GGa$, design parameters $\nu_1 = 1$, $\nu_2 =1.3$,
    $\alpha = 10$, $\beta = 40$ and $\eps = 0.0086$, and time-varying
    total load.  The example depicts the input-to-state stability of
    the mismatch dynamics.  }
  \label{fig:bursts}
\end{figure}

Next, we consider a time-varying total load given by a constant plus a
sinusoid, $P_l(t) = 4300 + 100 \sin(0.05t)$. With the same
communication topology, design parameters, and initial condition as
above, Figure~\ref{fig:one} (d)-(f) illustrates the behavior of the
network under the \dacplap dynamics. As established in
Proposition~\ref{pr:mismatch-ult-bound}, the total generation tracks
the time-varying load signal and the mismatch between these values has
an ultimate bound. Additionally, to illustrate how that the mismatch
vanishes if the load becomes constant, we show in
Figure~\ref{fig:bursts} a load signal that consists of short bursts of
sinusoidal variation that decay exponentially. The difference between
generation and load becomes smaller and smaller as the load tends
towards a constant signal.

Our final scenario considers addition and deletion of generators. The
initial communication topology is the undirected graph~$\GGau$
described in Table~\ref{tb:graphs}.  The design parameters and the
initial condition are the same as above.  The total load is $4200$ and
is same at all times. For the first $100$ seconds, the power
allocations converge to a neighborhood of a solution of the \ED
problem for the set of generators in $\GGau$. At time $t = 100s$, the
units $\{4, 11, 25, 45\}$ stop generating power and leave the
network. We select these generators because of their substantial
impact in the total power generation.  After this event, the resulting
communication graph is $\GGi$, cf.  Table~\ref{tb:graphs}.  The
generators implement the \tinv routine, after which the \dacplap
dynamics drives the mismatch to zero and minimizes the total cost. At
$t_2 = 200s$, another event occurs, the units $\{11, 45\}$ get added
to the network while the generator $27$ leaves. The resulting
communication topology is $\GGf$, cf.  Table~\ref{tb:graphs}. After
executing the \tinv routine, the dynamics converges eventually to the
optimizers of the \ED problem for the set of generators in $\GGf$, as
shown in Figure~\ref{fig:one}(g)-(i).  This example illustrates the
robustness of the \dacplap dynamics against intermittent generation by
the units, as formally established in
Proposition~\ref{pr:intermittent-gen}. In addition to the presented
examples, we also successfully simulated scenarios of the kind
described in Remark~\ref{re:load-agg-2}, where the total load is not
known to a single generator and is instead the aggregate of the local
loads connected to each of the generator buses, but we do not report
them here for space reasons.


\section{Conclusions}\label{sec:conclusions}

We have designed a novel provably-correct distributed strategy that
allows a group of generators to solve the economic dispatch problem
starting from any initial power allocation. Our algorithm design
combines elements from average consensus to dynamically estimate the
mismatch between generation and desired load and ideas from
distributed optimization to dynamically allocate the unit generation
levels.  Our analysis has shown that the mismatch dynamics between
total generation and load is input-to-state stable and, as a
consequence, the coordination algorithm is robust to initialization
errors, time-varying load signals, and intermittent power 
generation.  Our technical approach relies on tools from algebraic
graph theory, dynamic average consensus, set-valued dynamical systems,
and nonsmooth analysis, including a novel refinement of the LaSalle
Invariance Principle for differential inclusions that we have stated
and proved.  Future work will explore the study of the preservation of
the generator box constraints under the proposed coordination
strategy, the extension to scenarios that involve additional
constraints, such as transmission losses, transmission line capacity
constraints, ramp rate limits, prohibited operating zones, and
valve-point loading effects, and the study of the stability and
convergence properties of algorithm designs that combine our approach
here with traditional primary and secondary generator controllers.


{\scriptsize
\bibliographystyle{plainnat}%
\bibliography{biblio}
}

\appendix

\renewcommand{\theequation}{A.\arabic{equation}}
\renewcommand{\thetheorem}{A.\arabic{theorem}}

\section{Refined LaSalle Invariance Principle for differential
  inclusions}\label{app:refined-LaSalle}



In this section we provide a refinement of the LaSalle Invariance
Principle for differential inclusions, see e.g.,~\citep{JC:08-csm-yo},
by extending the results of~\citep{AA-CE:10} for differential
equations. Our motivation for developing this refinement comes from
the need to provide the necessary tools to tackle the convergence
analysis of the coordination algorithms presented in
Sections~\ref{sec:centralized}
and~\ref{sec:robust-dist-edp}. Nevertheless, the results stated here
are of independent interest.

\begin{proposition}\longthmtitle{Refined LaSalle Invariance Principle
    for differential inclusions}\label{pr:refined-lasalle-nonsmooth}
  Let $\setmap{F}{\real^n}{\real^n}$ be upper semicontinuous, taking
  nonempty, convex, and compact values at every point $x \in
  \real^n$. Consider the differential inclusion $\dot x \in F(x)$ and
  let $t \mapsto \varphi(t)$ be a bounded solution whose omega-limit
  set $\Omega(\varphi)$ is contained in $\SS \subset \real^n$, a
  closed embedded submanifold of $\real^n$.  Let $\OO$ be an open
  neighborhood of $\SS$ where a locally Lipschitz, regular function
  $\map{W}{\OO}{\real}$ is defined.  Assume the following holds,
  \begin{enumerate}
  \item the set $\EE = \setdef{x \in \SS}{0 \in \SetLie_F W(x) }$
    belongs to a level set of $W$,
    \label{as:refined-lasalle-1}
  \item for any compact set $\MM \subset \SS$ with $\MM \cap \EE =
    \emptyset$, there exists a compact neighborhood $\MM_c$ of $\MM$
    in $\real^n$ and $\delta < 0$ such that $\sup_{x \in \MM_c} \max
    \SetLie_{F} W(x) \le \delta$.
    \label{as:refined-lasalle-2}
  \end{enumerate}
  Then, $\Omega(\varphi) \subset \EE$.
\end{proposition}
\smallskip 

Before proceeding with the proof of the result, we establish an
auxiliary result.

\begin{lemma}\label{le:non-empty-intersection}
  Under the hypotheses of
  Proposition~\ref{pr:refined-lasalle-nonsmooth}, the sets
  $\Omega(\varphi)$ and $\EE$ have nonempty intersection.
\end{lemma}
\begin{pf}
  By contradiction, assume $\Omega(\varphi) \cap \EE =
  \emptyset$. Then, using the hypothesis~\ref{as:refined-lasalle-2} in
  Proposition~\ref{pr:refined-lasalle-nonsmooth}, there exists $\delta
  < 0$ such that $\sup_{x \in \Omega(\varphi)} \max \SetLie_F W(x) \le
  \delta$. Let $x \in \Omega(\varphi)$. Since this set is weakly
  positively invariant, there exists a trajectory $t \mapsto
  \tilde{\varphi}(t)$ of the differential inclusion with
  $\tilde{\varphi}(0) = x$ such that $\tilde{\varphi}(t) \in
  \Omega(\varphi)$.  Since $\frac{d}{dt} W(\tilde{\varphi}(t)) \in
  \SetLie_F W(\tilde{\varphi}(t))$ for almost all $t \ge 0$, we get
  $W(\tilde{\varphi}(t)) - W(x) \le \delta t$. This is in
  contradiction with the fact that $t \mapsto \tilde{\varphi}(t)$
  belongs to the compact set $\Omega(\varphi)$, where $W$ is lower
  bounded. \qed
\end{pf}

We are now ready to prove
Proposition~\ref{pr:refined-lasalle-nonsmooth}.

\begin{pf}[Proof of Proposition~\ref{pr:refined-lasalle-nonsmooth}]
  We consider two cases, depending on whether the set
  $\Omega(\varphi)$ (a) is or (b) is not contained in a level set
  of~$W$.  In case (a), given any $x \in \Omega(\varphi)$, there
  exists a trajectory of~$F$ starting at $x$ that remains in
  $\Omega(\varphi)$ (because of the weak positive invariance of the
  omega-limit set). If $x \not \in \EE$, then by the
  hypotheses~\ref{as:refined-lasalle-2}, there exists a compact
  neighborhood $\MM_x$ of $x$ in $\real^n$ and $\delta < 0$ such that
  $\sup_{y \in \MM_x} \SetLie_F W(y) \le \delta$. Since
  $\Omega(\varphi) \subset \SS$, the trajectory of $F$ starting at $x$
  remains in the set $\MM_x \cap \SS$ for a finite time, say $t_1$.
  Over the time interval $[0,t_1]$, we have $W(t) - W(0) \le \delta
  t$. This, however, is in contradiction with the fact that the
  trajectory belongs to $\Omega(\varphi)$ which is contained in a
  level set of $W$. Therefore, $x \in \EE$, and since this point is
  generic, we conclude $\Omega(\varphi) \subset \EE$.
  
  \begin{figure}
    \centering
    \includegraphics[width = 0.51
    \linewidth]{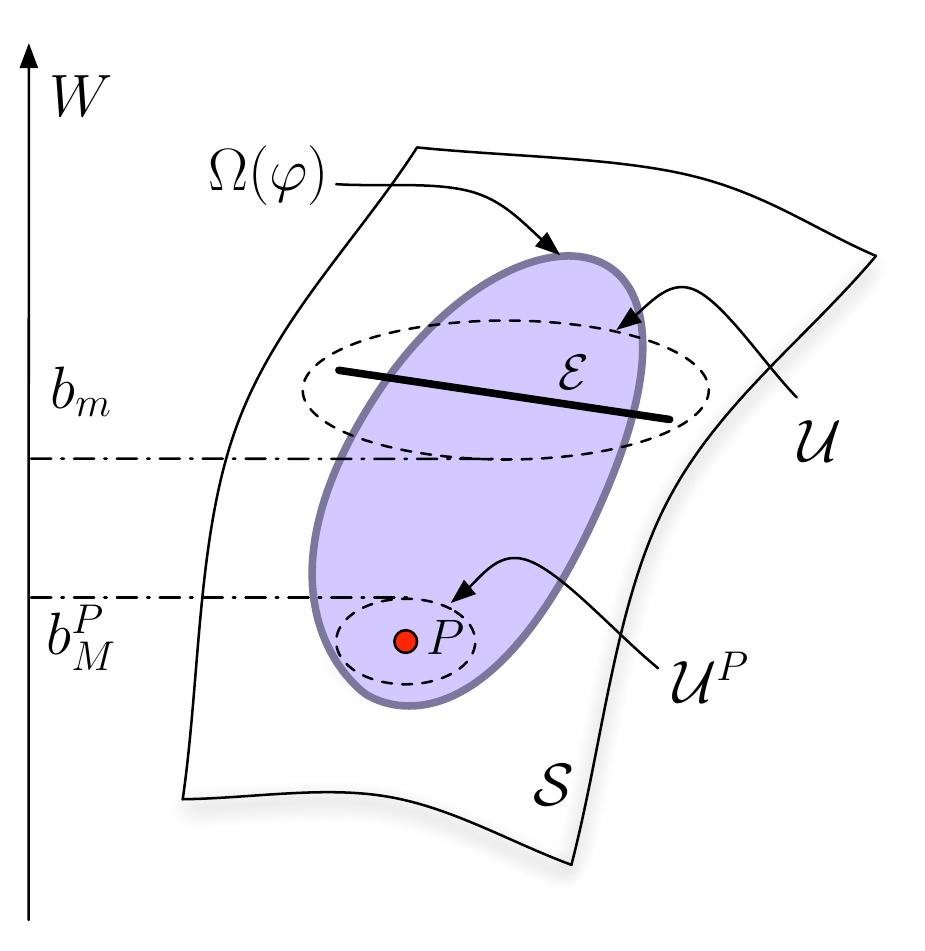}
    \includegraphics[width = 0.47
    \linewidth]{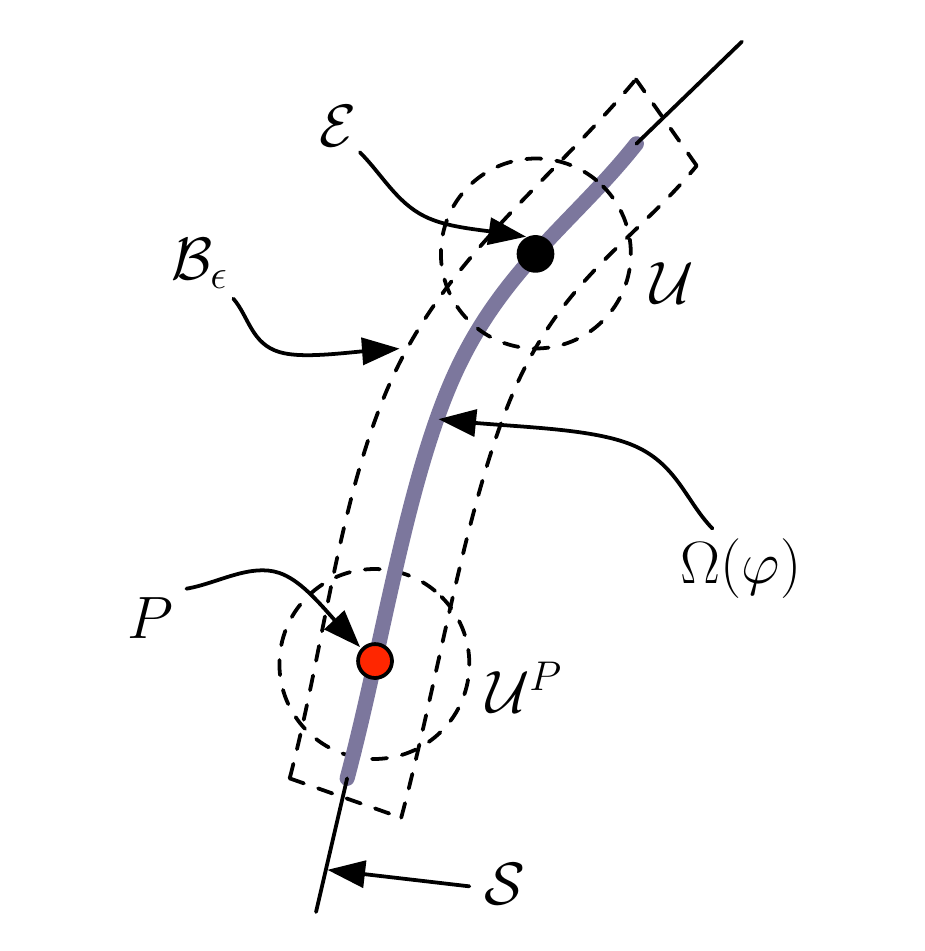}
    \caption{Illustration (adapted from~\citep[Figure 1]{AA-CE:10})
      depicting various elements involved in the case (b) of the proof
      of Proposition~\ref{pr:refined-lasalle-nonsmooth}.}\label{fig:refined}
  \end{figure}

  Next, we consider case (b) and reason by contradiction, i.e., assume
  that $\Omega(\varphi)$ is not contained in $\EE$ (see
  Figure~\ref{fig:refined}).  Given $\eps > 0$, let $\BB_\eps \subset
  \OO$ be a compact neighborhood of $\Omega(\varphi)$ in $\real^n$
  such that $d(\BB_\eps, \Omega(\varphi)) \le \eps$.  Let $\UU$ be an
  open neighborhood of $\EE$ in $\real^n$ and define $\UU_\eps = \UU
  \cap \BB_\eps$. Note that $\UU_\eps$ is nonempty because
  $\Omega(\varphi) \cap \EE$ is nonempty by
  Lemma~\ref{le:non-empty-intersection}.  Since $\Omega(\varphi)$ is
  not contained in a level set of $W$ but $\EE$ is by
  hypotheses~\ref{as:refined-lasalle-1}, we can choose $P \in
  \Omega(\varphi) \setminus \EE$ such that $W(P) \not =
  W(\EE)$. Without loss of generality, assume $W(P) < W(\EE)$ (the
  reasoning is analogous for the other case).  Select an open
  neighborhood $\UU^P$ of $P$ in $\real^n$ and define $\UU^P_\eps =
  \UU^P \cap \BB_\eps$. Define the following quantities
  \begin{alignat*}{2}
    b_m & = \inf_{x \in \UU_\eps} W(x) ,& \quad b^P_M & = \sup_{x
      \in \UU^P_\eps} W(x) .
  \end{alignat*}
  Note that the neighborhoods $\UU$ and $\UU^P$ can be chosen such
  that the set $\Omega(\varphi) \setminus (\UU \cup \UU^P)$ is
  nonempty, compact, and its intersection with $\EE$ is empty.  Along
  with this, one can select $\eps$ in such a way that $b_m > b^P_M$
  and from assumption~\ref{as:refined-lasalle-2} we get
  \begin{equation}\label{eq:Wdot-bound}
    \sup_{x \in \BB_\eps \setminus (\UU \cup \UU^P)}
    \max \SetLie_{F} W(x) \le \delta <0 ,
  \end{equation}
  (in the case $W(P) > W(\EE)$, we would reason with the quantities
  $b_M = \sup_{x \in \UU_\eps} W(x)$ and $b_m^P = \inf_{x \in
    \UU_\eps^P} W(x)$).
  Since $\Omega(\varphi)$ is the omega-limit set of $\varphi$ and
  $\BB_\eps$ is a compact neighborhood of $\Omega(\varphi)$, there
  exists $t_1 > 0$ such that $\varphi(t_1) \in \UU^P_\eps$ and
  $\varphi(t) \in \BB_\eps$ for all $t \ge t_1$.  Moreover, since
  $\Omega(\varphi) \cap \EE$ is nonempty, there must also exist $t_2 >
  t_1$ such that $\varphi(t_2) \in \UU_\eps$. From continuity of the
  trajectory we deduce that there exist times $t_1^*, t_2^* \in (t_1,
  t_2)$, $t_1^* < t_2^*$ such that $\varphi(t_1^*)$ and
  $\varphi(t_2^*)$ lie on the boundary of the compact set $\BB_{\eps}
  \setminus (\UU_{\eps} \cup \UU_{\eps}^P)$, with $\varphi(t_1^*)$
  belonging to the closure of $\UU^P_\eps$ and $\varphi(t_2^*)$ to the
  closure of $\UU_\eps$.  However, this is not possible as
  $W(\varphi(t_2^*)) \ge b_m > b^P_M \ge W(\varphi(t_1^*))$ and, in
  the interval $[t_1^*,t_2^*]$, the trajectory belongs to $\BB_{\eps}
  \setminus (\UU_{\eps} \cup \UU_{\eps}^P)$, where the function~$W$
  can only decrease due to~\eqref{eq:Wdot-bound}, which is a
  contradiction. \qed
\end{pf}

\section{Continuity properties of set-valued Lie
  derivatives}\label{app:continuity-Lie}

Here we present an auxiliary result for the convergence analysis of
the algorithms of Sections~\ref{sec:centralized}
and~\ref{sec:robust-dist-edp}.

\begin{lemma}\longthmtitle{Continuity property of set-valued Lie
    derivatives}\label{le:continuity-set-valued-gen}
  Let $\map{W}{\real^n}{\real}$ be a locally Lipschitz and regular
  function.  Let $\map{g}{\real^n \times \real^n}{\real^n}$ be a
  continuous function and define the set-valued map
  $\setmap{F}{\real^n}{\real^n}$ by $F(x) = \setdef{g(x,\zeta)}{\zeta
    \in \partial W(x)}$.  Assume that
  \begin{enumerate}
  \item $\SS$ is an embedded submanifold of $\real^n$ such that
    $\zeta^\top g(x,\zeta) \le 0$ for all $x \in \SS$ and all $\zeta
    \in \partial W(x)$, \label{as:continuity-1}
  \item for any $x \in \SS$, if $\zeta^\top g(x,\zeta) = 0$ for some
    $\zeta \in \partial W(x)$, then $x \in \EE = \setdef{z \in \SS}{0
    \in \SetLie_F W(z)}$. \label{as:continuity-2}
  \end{enumerate}
  Then, for any compact set $\MM \subset \SS$ with $\MM \cap \EE =
  \emptyset$, there exists a compact neighborhood $\MM_c$ of $\MM$ in
  $\real^n$ and $\delta < 0$ such that $ \sup_{x \in \MM_c} \max
  \SetLie_F W(x) \le \delta$.
\end{lemma}
\begin{pf}
  We reason by contradiction, i.e., assume that for all compact
  neighborhoods $\MM_c$ of $\MM$ in $\real^n$ and all $\delta < 0$, we
  have
  \begin{align*}
    \sup_{x \in \MM_c} \max \SetLie_{F} W(x) > \delta .
  \end{align*}
  Note that this implies that $ \sup_{x \in \MM_c} \max \SetLie_{F}
  W(x) \ge 0$.  Now, for each $k \in \integerspositive$, consider the
  compact neighborhood $\MM_k = \MM + \overline{B(0,\frac{1}{k})}$
  of~$\MM$. From the above, we deduce the existence of a sequence
  $\{x_{k}\}_{k=1}^{\infty}$ with $x_{k} \in \MM_k$ such that
  \begin{equation}\label{eq:max-limit3}
    \lim_{k \to \infty} \max \SetLie_{F}
    W(x_{k}) = \ell \ge  0.
  \end{equation}
  Since the whole sequence belongs to the compact set $\MM_1$, there
  exists a subsequence, which we denote with the same indices for
  simplicity, such that
  \begin{equation}\label{eq:point-limit3}
    \lim_{k \to \infty} x_{k} = \tilde{x} \in \MM .
  \end{equation}
  From~\eqref{eq:max-limit3}, there exists a sequence $\zeta_{k}
  \in \partial W(x_{k})$ such that
  \begin{align}\label{eq:set-lie-limit}
    \lim_{k \to \infty} \zeta_{k}^\top g(x_{k},\zeta_{k}) \ge 0.
  \end{align}
  Since $\partial W$ is upper semicontinuous with compact values, the
  set $\partial W(\MM_1)$ is compact, cf.~\citep[Proposition 3,
  p. 42]{JPA-AC:84}.  This implies that the sequence $\{\zeta_{k}\}$
  belongs to the compact set $\partial W(\MM_1)$ and so, there
  exists a subsequence, denoted again by the same indices for
  simplicity, such that $\zeta_{k} \to \tilde{\zeta}$. Since
  $\partial W$ is upper semicontinuous and takes closed values, we
  deduce from~\citep[Proposition 2, p. 41] {JPA-AC:84} that
  $\tilde{\zeta} \in \partial W(\tilde{x})$.
  From~\eqref{eq:point-limit3} and~\eqref{eq:set-lie-limit}, since $g$
  is continuous, we obtain $\tilde{\zeta}^\top
  g(\tilde{x},\tilde{\zeta}) \ge 0$.  By
  assumption~\ref{as:continuity-1}, this implies
  $\tilde{\zeta}^\top g(\tilde{x},\tilde{\zeta}) = 0$.
  Assumption~\ref{as:continuity-2}
  then implies $\tilde{x} \in \EE$, which along
  with~\eqref{eq:point-limit3} contradicts~$\MM \cap \EE = \emptyset$.
  \qed
\end{pf}

\end{document}